\documentclass[12pt]{article} 
\usepackage[dvips]{graphics,color}
\usepackage{latexsym}
\usepackage{amssymb}

\let\cal\mathcal

\usepackage{latexsym}
\usepackage{amssymb}
\usepackage{euscript}
\usepackage[dvips]{graphics}
\usepackage{epsf}
\let\cal=\mathcal      

\def\mcc{M\raise.5ex\hbox{c}C}
\def\mccarthy{M\raise.5ex\hbox{c}Carthy}

\def\eg{{\it e.g. }}
\def\ie{{\it i.e. }}

\def\h{{\cal H}}

\def\K{{\cal K}}
\def\M{{\cal M}}

\def\hid{H^\infty(\D)}

\def\a{\alpha}
\def\l{\lambda}
\def\z{\zeta}

\def\vare{\varepsilon}

\let\i=\infty

\def\={\ = \ }

\def\A{{\cal A}}

\def\CC{{\cal C}}

\def\C{\mathbb C}
\def\R{\mathbb R}
\def\T{\mathbb T}
\def\D{\mathbb D}

\def\Z{\mathbb Z}
\def\NN{\mathbb N}

\def\inn{\ \in \ }

\def\be{\setcounter{equation}{\value{theorem}} \begin{equation}}
\def\ee{\end{equation} \addtocounter{theorem}{1}}
\def\beq{\begin{eqnarray*}}
\def\eeq{\end{eqnarray*}}
\def\se{\setcounter{equation}{\value{theorem}}} 
\def\att{\addtocounter{theorem}{1}}
\def\vs{\vskip 5pt}
\def\bs{\vskip 12pt}

\def\rem{\bs \att {{\bf Remark \thetheorem \ }} }
\def\bp{{\sc Proof: }}
\def\ep{{}{\hfill $\Box$} \vskip 5pt \par}

\def\bl{\begin{lemma}}
\def\el{\end{lemma}}
\def\bt{\begin{theorem}}
\def\et{\end{theorem}}
\def\bprop{\begin{prop}}
\def\eprop{\end{prop}}
\def\bd{\begin{definition}}
\def\ed{\end{definition}}
\def\br{\begin{remark}}
\def\er{\end{remark}}
\def\bexer{\begin{exercise}}
\def\eexer{\end{exercise}}
\def\bfig{\begin{figure}}
\def\efig{\end{figure}}

\newtheorem{theorem}{Theorem}[section]
\newtheorem{prop}[theorem]{Proposition}
\newtheorem{lemma}[theorem]{Lemma}
\newtheorem{cor}[theorem]{Corollary}

\newtheorem{question}[theorem]{Question}
\newtheorem{definition}[theorem]{Definition}

\newcommand{\U}{\cal U}

\newcommand{\al}{\alpha}
\newcommand{\hac}{hyperbolic\ algebraic\ curve\ }
\newcommand{\hacp}{hyperbolic\ algebraic\ curve}
\newcommand{\hacs}{hyperbolic\ algebraic\ curves\ }
\newcommand{\hacsp}{hyperbolic\ algebraic\ curves}
\newcommand{\gasp}{hyperbolic\ analytic\ curve}
\newcommand{\gass}{hyperbolic\ analytic\ curves\ }
\newcommand{\gassp}{hyperbolic\ analytic\ curves}
\newcommand{\gas}{hyperbolic\ analytic\ curve\ }
\newcommand{\cn}{\C^n}
\newcommand{\ol}{1/\overline{\l}}
\newcommand{\om}{1/\overline{\mu}}
\newtheorem{example}[theorem]{Example}
\newcommand{\EE}{\mathbb E}
\newcommand{\PPP}{\mathbb P}
\newcommand{\alg}{{\cal Alg}}

\newcommand{\CP}{{\C \PPP}}
\newcommand{\Dp}{\D^+}
\newcommand{\Dm}{\D^-}
\newcommand{\zp}{0^+}
\newcommand{\zm}{0^-}
\renewcommand{\sp}{\sigma^+}
\newcommand{\sm}{\sigma^-}
\newcommand{\precc}{\prec\ }
\newcommand{\Mp}{\M_{+}}
\newcommand{\Mn}{\M_{-}}
\newcommand{\Ep}{E_{+}}
\newcommand{\Em}{E_{-}}
\newcommand{\Wp}{W_{+}}
\newcommand{\Wm}{W_{-}}

\renewcommand{\P}{{\mathcal P}}

\begin{document}
\setlength{\baselineskip}{21pt}
\title{ Hyperbolic Algebraic and Analytic Curves}
\author{Jim Agler
\thanks{Partially supported by National Science Foundation Grant
DMS 0400826}\\
U.C. San Diego\\
La Jolla, California 92093
\and
John E. M\raise.5ex\hbox{c}Carthy
\thanks{Partially supported by National Science Foundation Grant
DMS 0501079}\\
Washington University\\
St. Louis, Missouri 63130}
%\date{October 27 2004}

\bibliographystyle{plain}

\maketitle
\begin{abstract}
A \hac is a bounded subset of an algebraic set. We study the function theory 
and functional analytic aspects of these sets. We show that their function theory
can be described by  finite codimensional subalgebras of the holomorphic functions on 
the desingularization. We show that classical analytic techniques, such as interpolation, 
can be used to answer geometric questions about the existence of biholomorphic maps.
Conversely, we show that the algebraic-geometric viewpoint leads to interesting 
questions in classical analysis.
\end{abstract}

\baselineskip = 18pt

\setcounter{section}{-1}
\section{Introduction}\label{seca}
By a {\em \hac} we mean 
%an open subset $V$ of a one dimensional algebraic set $\CC$ in
%$\cn$ that has a piecewise smooth boundary and such that every component of $V$ supports
%non-constant bounded holomorphic functions. The obvious way to attain
%such a set is to take $V$ to be the intersection of $\CC$ with some
%semi-algebraic bounded open set $U$ in $\cn$, and these are the examples in which we
%are primarily interested. 
a set $V$ that is the intersection of a  one dimensional algebraic set $\CC$ 
with a bounded open set in $\C^n$.
We call them hyperbolic to emphasize that 
we wish to study the geometry and function theory on $V$, not the global
theory on $\CC$, nor just the local theory on small neighborhoods of
points of $V$. Note that by a holomorphic function on $V$ we mean a
function that in a neighborhood of every point of $V$ is the restriction
of a holomorphic function on a neighborhood in $\cn$.

The leimotif is the difference between studying the hyperbolic geometry
and function theory on the unit disk $\D$ in $\C$, as opposed to
studying the ``parabolic'' theory on $\C$ or the ``spherical'' theory on 
the Riemann sphere. Here are some simple examples of what we wish to consider.

%\exam \label{exa1}
% {\bf Example 1}

\begin{example}\label{exa1}
{\rm
Let $\CC$ be Newton's nodal cubic, $\{ (z,w) \inn \C^2 \ : \
z^2 = w^2 (1-w)\}$, and let $V_1 = \CC \cap \{ |1-w| < 2 \}$ (see
Figure~\ref{pica1}).
\bfig 
\resizebox{!}{2in}{\includegraphics {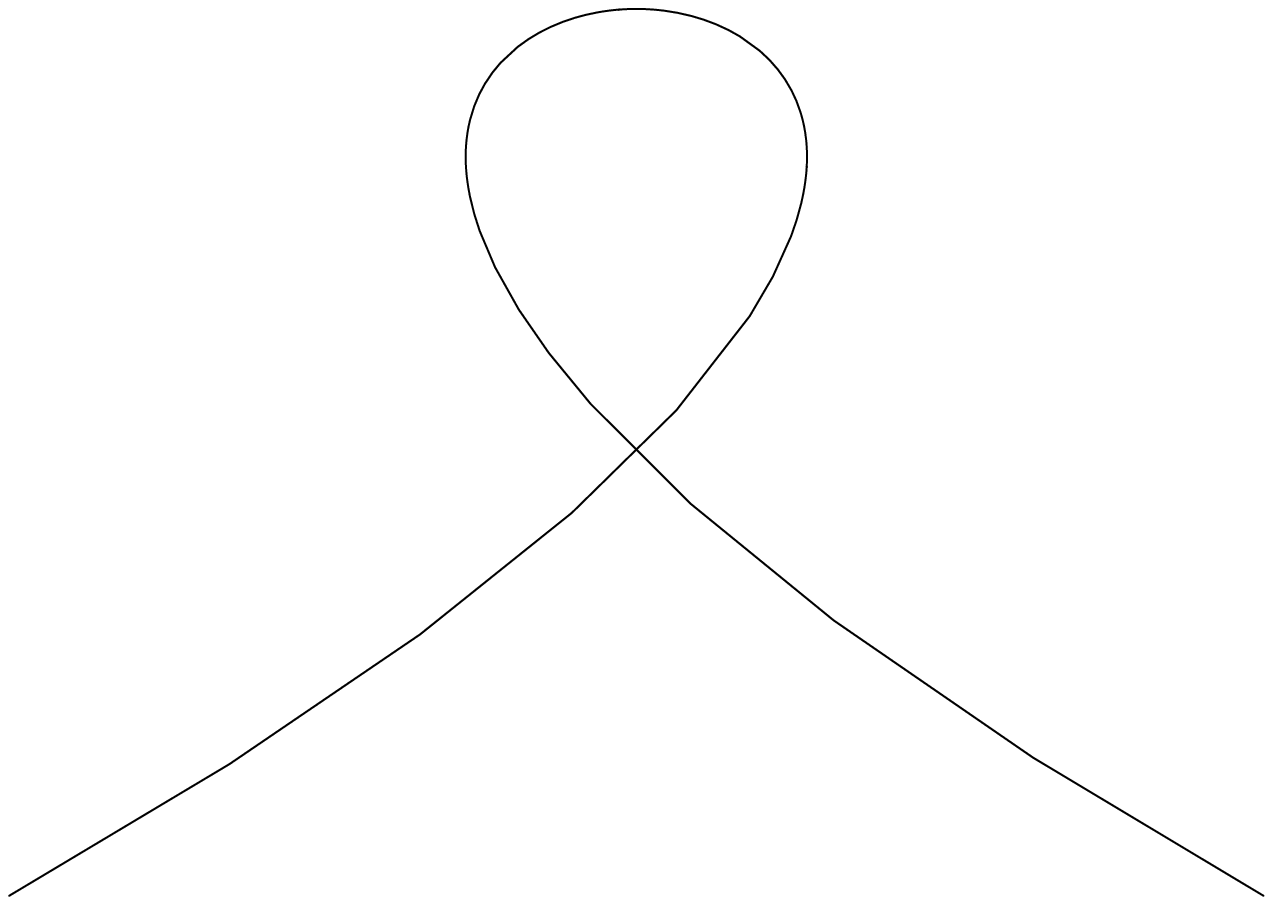}}
\caption{$V_1$, a hyperbolic subset of Newton's nodal cubic} 
\label{pica1} \efig
Now consider another \hac $V_2$ that is  topologically
the same as $V_1$, \ie a disk that intersects itself once. 
When is there a biholomorphic bijection from $V_1$
onto $V_2$? 

The answer (see Section~\ref{secd}) is that the Kobayashi geodesic that goes from
the multiple point across $V_2$ and back to the multiple point must have
the same length as the corresponding geodesic in $V_1$. If one cuts
$\CC$
very close to the multiple point, this length gets very
large; if one cuts it very far away, the length is small. So there is a
one parameter family of conformal structures on sets that are
homeomorphic to $V_1$, and they can all be attained by taking the
intersection  of $\CC$ with an appropriate bounded set $\U$. 

The set $V_1$ 
%in Example~\ref{exa1} 
can be obtained as the image of the
unit disk $\D$ under the map
$$
h(\z) \= (\sqrt{2}(\z - 2\z^3), 1 - 2\z^2).
$$
The map $h$ has nowhere vanishing derivative, and it is one-to-one except 
at the points $\pm 1/\sqrt{2}$.

The set $\{ F \circ h \ : \ F\ {\rm holomorphic\ on\ } V_1 \}$
is the subalgebra of $O(\D)$ (the functions analytic on $\D$) that have the same value at 
$\pm 1/\sqrt{2}$. One could reverse the viewpoint, and ask:

Given a finite codimensional subalgebra $A$ of 
$O(\D)$, when does there exist a map $h$ from $\D$ to some \hac $V$ in $\C^n$
such that
\be
\label{eqa01}
A \= \{ F \circ h \ : \ F\ {\rm holomorphic\ on\ } V \} ?
\ee

If such an $h$ exists, what is the smallest dimension $n$ that suffices?
}
\end{example}
\vs
%\exam \label{exa2}
%{\bf Example 2}
\begin{example}\label{exa2}
{\rm
%Let $V$ be three disks in $\C^3$ that intersect at a common point.
%When is $V$ biholomorphically equivalent to a \hac in $\C^2$?
%
%A moment's thought answers the question. If one moves the common point
%to the origin, and writes the three disks as
%$$
%D_i \= \{ z \l_i \ : \ z \inn \D \},\qquad i = 1,2,3,
%$$
%for $\l_i$ non-zero vectors in $\C^3$, then $V$ can be embedded in
%$\C^2$ if and only if $\{ \l_1,\l_2,\l_3 \}$ are linearly dependent.
Let $A^2$ and $A^3$ be the subalgebras of $O(D)$ given by the following
(the superscript is the codimension):
\beq
A^2 &\=& \{ f \inn O(\D) \, : \, f(0) = f(\frac{1}{2}) = f(-\frac{1}{2}) \} \\
A^3 &\=& \{ f \inn A^2 \, : \, f'(0) \= -\frac{15}{64} \left[ f'(\frac{1}{2}) + f'(-\frac{1}{2}) \right] \}.
\eeq
Can these algebras be realized as in (\ref{eqa01}) for some $V$ in $\C^2$?

The algebra $A^2$ cannot; because of the triple point, for any $h$ mapping into $\C^2$, there must be 
some linear relation between $h(0), h(\frac{1}{2})$ and $h(-\frac{1}{2})$, and this relation will
apply to the entire algebra (\ref{eqa01}) generated by $h$. One can realize $A^2$ in $\C^3$ --- 
see Section~\ref{sece}.

Once one adds a linear relation on the derivatives, as in $A^3$, one can realize the algebra in $\C^2$.
For the particular example of $A^3$, the function
$$
h(z) \= (z \frac{z^2 - 1/4}{1 - z^2/4},\ z^2 \frac{z^2 - 1/4}{1 - z^2/4} )
$$
works.
}
\end{example}

\vs
Examples~\ref{exa1} 
 and 
\ref{exa2} 
 illustrate the first three of the general themes of
this paper: 

(i) deciding when two \hacs are isomorphic (in the sense of
there
being a biholomorphic homeomorphism between them);

(ii) studying the duality between \hacs and cofinite subalgebras
of the holomorphic functions on the desingularization;

(iii) 
finding the minimum dimension into which a given \hac can be
isomorphically realised;

(iv)
exploring the operator theoretic aspects of \hacsp.
\vs
%The set $V_1$ in Example~\ref{exa1} 
%can be obtained as the image of the
%unit disk $\D$ under the map
%$$
%h(\z) \= (r\z - r^3\z^3, 1 - r^2\z^2),
%$$
%where $r =\sqrt{2}$.
%
%The map $h$ has nowhere vanishing derivative, and it is one-to-one except 
%at the points $\pm 1/\sqrt{2}$.

There is in general a duality between \hacs and maps like $h$ that are
one-to-one except on  a finite set.
To describe the duality formally, it
is convenient to broaden our scope to the
analytic analogue of a \hacp, which we call a {\em \gasp.}
It is a one dimensional  analytic subvariety of a bounded open set.

\bd
\label{defa1}
A \gas $V \subseteq \C^n$ is a set such that

\noindent
(i) $V$ is a relatively closed subset of
 some bounded open set $U \subseteq \C^n$.

\noindent
(ii) There is a finite set $B \subset V$ such that
for every $P\ \in \ V\setminus B$ there exists a neighborhood $W$ of $P$
in $\C^n$ and $n-1$ holomorphic functions $F = (F_2, \dots, F_n)$
on $W$ such that 

(a) The rank of the derivative $D F$ is $n-1$ on $W$.

(b) $V \cap W \= \{ z \ \in \ W\ : \ F(z) = 0 \}$.

\noindent
(iii) For every $P\ \in \ B$ there exists a neighborhood
$W$ of $P$
in $\C^n$ and finitely many   holomorphic tuples
of functions $F^j = (F^j_2, \dots, F^j_{n_j})$
on $W$ such that
$$
V \cap W \= \cup_j \{ z \ \in \ W\ : \ F^j(z) = 0 \} .$$
\ed

\bd
\label{defa2}
A holomap is a proper holomorphic map $h$ from a 
Riemann surface $S$ 
into a bounded open set $U \subseteq\C^n$ such that
there is a finite subset $E$ of $S$ with the property that
 $h$ is non-singular
and injective on $S\setminus E$.
% and such that $h$ is non-constant on every component of $S$.
\ed
%(Note that if $S$ is hyperbolic,
%the restriction that $h$ be proper forces it to be
%non-constant on every component.)
% (as $S$ has no compact components).
\vs
%We show in Theorem~\ref{thmb1} that 
All \gass arise as the images of
Riemann surfaces under holomaps, and the surface is unique if it is
unpunctured (see Section~\ref{secb} for a precise definition of unpunctured).

\vs
%\begin{othm}~\ref{thmb1}
{\bf Theorem~\ref{thmb1}}{\em
\
(i) If $h$ is a holomap from a
 Riemann surface $S$ into the bounded open set $ U \subset \C^n$, then
$h(S)$ is a \gas.

(ii) Conversely,
if $V$ is a \gas relatively closed
in the bounded open set $U$,
then there exists a
 Riemann surface
$S$ and a 
holomap $h : S \to U \subset \subset \C^n$ such that $V = h(S)$.

(iii) Moreover, if there is another
Riemann surface $\tilde S$ and
a holomap $\tilde h$ from $\tilde S$ onto $V$, and if $\tilde S$
is unpunctured on the pre-image of an open set containing all the
singular points of $V$, then there is a
biholomorphic bijection $\psi: S \to \tilde S$ such that $ h =
\tilde h \circ \psi$.
}
%\end{othm}
\vs

One can also approach \gass algebraically.
If $h : S \to V$ is a holomap, then any holomorphic function $\phi$ on
$V$ can be lifted to the function $\phi \circ h$ on $S$. 
Not all functions on $S$ arise in this way: for example, such functions
must identify all points that $h$ identifies. But the set of functions
on $S$ that do factor through $V$ via $h$ is an algebra, of finite
codimension. Let us make two definitions.
\bd
If $W$ is any subset of a set $S$ with a conformal structure, then by 
a holomorphic function on $W$ we mean a function $f:W \to \C$ with the 
property that for every point $\l$ in $W$ there is a neighborhood $U$ of $\l$ in $S$ and 
a holomorphic function $\phi$ on $U$ such that $ \phi |_{U \cap S} = f|_{U \cap S}$.
We denote by $O(W)$ the Fr\'echet space of all holomorphic functions on $W$, 
with the topology of uniform convergence on compacta.
\ed
\bd
A {\em cofinite subalgebra} $A$ of $O(W)$ is a unital closed subalgebra of $O(W)$ 
whose codimension
as a vector space is finite.
When the set $W$ is understood, we shall call $A$ a cofinite algebra. 
\ed
We prove in Theorem~\ref{thmb3} that all cofinite subalgebras of $O(S)$,
where $S$ is good Riemann surface, arise as $\{ F \circ h  : F \inn O(V)\}$
for a holomap $h: S \to V$, and conversely that every set of
the form $\{ F \circ h : F \inn O(V)\}$
is a cofinite subalgebra. Moreover, if $S$ is a reasonably nice \hacp,
then $V$ can also be chosen to be a \hac.
So all cofinite algebras of $O(S)$ arise as in Example~\ref{exa1}
by factoring through another \hac that has a finite number of singularities.
\vs

 We shall
write
$$
A_h \ := \ \{ F \circ h \ : \ F \ \in \ O(V) \}.
$$
%Thus one can study \gassp, and \hacsp, algebraically too.

To state the theorem, we need one more definition.
\bd
A Riemann surface $S$ is called {\em good} if there is a holomap
$h:  S \to U \subseteq \C^n$ such that $h$ has no singular points and
$\{ p \circ h \ : p \, \in \, \C [z_1,\dots,z_n] \}$ is dense in $O(S)$.
\ed
\rem
Every finite Riemann surface is good, where a {\em finite Riemann
surface} is one obtained from  a compact Riemann surface by deleting a
positive 
finite number of smoothly bounded disjoint closed disks. Indeed, by a
result of E.L. Stout
\cite{sto66}, every finite Riemann surface $S$ can be properly embedded
in $\D^3$ by a triple of inner functions. The image $h(S)$
is an analytic subvariety of $\D^3$ by Theorem~\ref{thmb1},
and has no singular points, so every analytic function on $S$ is of
the form $F \circ h$ for some analytic function $F$ on a neighborhood
of $h(S)$. By a theorem of H.~Cartan \cite{car51},
$F$ can be extended to an analytic function on the whole tridisk,
and so can be approximated u.c.c. by polynomials.
Therefore $A_h = O(S)$.

\vs
We can now state the equivalence between cofinite algebras and 
realizations of them via holomaps.

{\bf Theorem~\ref{thmb3}}{\em
\
(i) If $h : S \to V \subseteq U$ is a holomap from an unpunctured
 Riemann
surface onto  a \gas $V$, then
$$
A_h \ := \ \{ F \circ h \ : \ F \ \in \ O(V) \}
$$
is a cofinite subalgebra of $O(S)$.

(ii) Conversely, if $A$ is a cofinite subalgebra of $O(S)$ for some
good Riemann surface $S$, then there
exists a holomap $h : S \to \C^n$ such that $A = A_h$.

(iii) Let  $V$ be a \hacp, and assume that $O(V)$ can be generated by a 
finite number of rational functions that are continuous on $\overline{V}$.   
If $A$ is a cofinite subalgebra of $O(V)$,
then there is a \hac $W$ and a holomap $h$ from $V$ onto $W$ such that
$A = A_h$.
}
\vs
The purpose of proving these equivalences is to get the isomorphism
theorem~\ref{thmcc1}, which, roughly speaking, says that two \hacs (or \gassp) are
isomorphic (in the sense of there being a biholomorphic homeomorphism between them)
if and only if there is an isomorphism of the 
corresponding Riemann surface that preserves the structure of the cofinite algebras.
\vs
If $A = A_h$ for some holomap $h$, we shall say that $h$ {\em holizes}
the algebra $A$, and that it holizes the set $h(S)$.
\bs
In Sections~\ref{secy}-\ref{sece}, we consider cofinite subalgebras of $O(\D)$, which we
call petals. 
It is of particular interest to know when cofinite subalgebras of $O(\D)$ can be holized by
a pair of finite Blaschke products. The \hac is then a special sort of subvariety $V$ of the bidisk that 
has the property that $\overline{V} \cap \partial (\D^2) \= \overline{V} \cap (\partial \D)^2$.
These are called distinguished varieties (since they exit the disk through the distinguished boundary)
and have been studied in 
\cite{agmc_dv} and \cite{veg07}.

One difficulty in the theory is proving that a particular pair of functions in a cofinite algebra 
generate the whole algebra. Thanks to work of J.~Wermer \cite{wer58}, petals are easier to understand than
algebras on general surfaces. We explain this in Section~\ref{secy}.
For a pair of Blaschke products, there is a convenient formula that gives the exact
codimension of the algebra they generate.

{\bf Theorem~\ref{thmbp2}}{\em
\
Let $f$ and $g$ be Blaschke products of degree $m$ and $n$ respectively.
Suppose there are exactly $r < \infty$ unordered pairs of points on $\T$ that are not separated by $(f,g)$.
Then the codimension of $\alg(f,g)$ is
$\frac{(m-1)(n-1) - r}{2} $.
}

In Section~\ref{secz}, we show how the question of whether there is a holomorphic map from one
petal to another is reducible to a Pick interpolation problem.
In Section~\ref{secd}, we describe all petals of codimension one. In Section~\ref{sece}, we discuss petals of codimension
two, and our results here are not complete.

A detailed study of the bounded analytic functions on the unit disk, the space $\hid$,
and the Hardy space $H^2$ has had profound consequences in operator theory.
The backward shift operator $S^\ast$ on $H^2$ is in some sense universal for contractions
(operators of norm less than or equal to $1$). A celebrated theorem of B.~Sz.-Nagy asserts
that every completely non-unitary contraction
can be realised as the restriction of a direct sum of copies of  $S^\ast$ to an invariant subspace \cite{szn53}.
See \eg \cite{szn-foi}, \cite{nik},  or \cite{ampi} for expositions.

In Section~\ref{seci}, we illustrate how every \hac comes with its own 
model theory. 
We consider $V \= \{ (z,w) \inn \D^2 \, : \, z^2 = w^2 \}$. We prove the following 
analogue of the von Neumann-Wold theorem.
{\bf Theorem~\ref{thmi2}}{\em
\
Let $T = (T_1,T_2)$ be a pair of commuting isometries on a Hilbert space $\h$ satisfying $T_1^2 = T_2^2$.
Then there is a decomposition 
$\h \= \M_{+} \oplus \M_{-} \oplus \K$, isometries $W_{\pm}$ on 
$\M_{\pm}$, and a pair of operators $E_{\pm} : \K \to \M_{\pm}$ satisfying
$\Ep^\ast \Ep + \Em^\ast \Em = I$, $\Wp^\ast \Ep = 0$, and $\Wm^\ast \Em = 0$,
so that
\beq
\nonumber
T_1 &\ \cong \ & 
\bordermatrix{&\Mp  &\Mn& \K\cr
\Mp &\Wp & 0&\Ep\cr
\Mn&0&\Wm &\Em\cr
\K&0&0&0} \\
&&\\
&&\\
T_2 &\ \cong \ & 
\bordermatrix{&\Mp  &\Mn& \K\cr
\Mp &\Wp & 0&\Ep\cr
\Mn&0&-\Wm &-\Em\cr
\K&0&0&0} .
\eeq
}

In Section~\ref{secnp}, we show how the geometry and analysis interact in finding the minimal
inner functions on a \hacp. We show that on the Neil parabola, the minimal inner functions split into
two sorts, one of degree 2 and the other of degree 3.

\section{Lifting \gassp}
\label{secb}

The following theorem is well-known, but we could
not find a convenient proof in the literature, so we include one for
completeness.
\bd
A puncture in a Riemann surface $S$ is a domain $D_0 \subseteq S$
that is conformally equivalent via a map $\tau$ to $\{ z \, \in \, \C \ : \
0 < |z| < 1 \}$, and such that every sequence $\{z_n \}$ with $\tau(z_n)
\to 0$ is discrete in $S$.
We say that the surface is unpunctured if it has no punctures.
\ed
\bt
\label{thmb1}
(i) If $h$ is a holomap from a 
 Riemann surface $S$ into the bounded open set $ U \subset \C^n$, then
$h(S)$ is a \gas.

(ii) Conversely,
if $V$ is a \gas in the bounded open set $U$,
then there exists a 
 Riemann surface
$S$ and a holomap $h : S \to U \subset \subset \C^n$ such that $V = h(S)$.

(iii) Moreover, if there is another
Riemann surface $\tilde S$ and
a holomap $\tilde h$ from $\tilde S$ onto $V$, and if $\tilde S$
is unpunctured on the pre-image of an open set containing all the
singular points of $V$, then there is a
biholomorphic bijection $\psi: S \to \tilde S$ such that $ h =
\tilde h \circ \psi$.
\et

\bp
(i) 
Because $h$ is proper, $h(S)$ must be relatively closed in $U$.
Let $V = h(S)$ and $B = h(E)$. For every point  $P$ in $V \setminus
B$,
the inverse image $Q = h^{-1}(P)$
is unique, and, after a change of basis if necessary,
we can find a local coordinate $\z$ at $Q$ so that locally
$$
h(\zeta) \= (\zeta, h_2(\z),\dots, h_n(\z) ).
$$
For $ 2 \leq i \leq n$, define
$F_i(z) = z_i - h_i(z_1) $. Then the derivative of $F$ has full
rank.

Now let $P \in B$. Let $Q$ be one of the (finitely many) pre-images
of $P$ under $h$. We will show that for some neighborhood $D$ of
$Q$, the set $h(D)$ is an analytic subvariety of $U$.

Indeed, after a change of basis if necessary, we can assume
that there is some local coordinate $\z$ at $Q$ and some positive
integer $N$ such that, on some disk $\Delta = \D(Q,\vare)$, the
function $h$ is injective and of the form
\be
\label{eq5}
h(\zeta) \= (\z^N, h_2(\z),\dots, h_n(\z) ).
\ee
(See \cite[p.125]{chi90}.) 
Consider the function $g: h(\Delta) \to D := \D(0,\vare^N)$ given by
projection onto the first coordinate. Then $g$ is an analytic
branched cover of branching order $N$ over $D$. (This means that
$g$ is proper, and that $g$ is an $N$-sheeted covering map of $D
\setminus \{ 0\}$). Therefore, by \cite[Thm. III.B.19]{gur},
$h(\Delta)$ is a one-dimensional 
analytic subvariety of a neigborhood of $P$,
which means it satisfies 
condition (iii) of Definition~\ref{defa1}. Taking the union over
all pre-images of $P$, we get that $h(S)$ is an analytic variety at
$P$ too.

\vs
(ii) Suppose $V$ is a \gas.
Let $B$ be the set of singular points of $V$.
Let $
S' = V \setminus B$
with the conformal structure inherited from $V$.
For each $P$ in $B$, $V$ is locally the union of finitely many
irreducible analytic varieties. 
Each of these irreducible pieces is locally the image of a disk
under a Puiseux series as in (\ref{eq5}).

%By the Weierstrass preparation theorem,
%\cite[p. 107]{fis01}
%for every point $P$ in $B$, there are topological disks $D_1, \dots
%$D_n$
%in $V$ such that $D_i \cap D_j = \{P \} $ for every $ i \neq j$.
%
%By assumption 
%(iv), 
%if $P$ is in $B$ then locally $S$ is
%$$
%S \= \{ z \ : \ F_i(z_1,z_i) = 0,\ 2\leq i \leq n \}.
%$$
%Factor each $F_i$ into irreducible Weierstrass polynomials,
%$$
%F_i \= {\rm (unit)} * \prod_{j_i=1}^{t_i} G_{i,j_i} .
%$$
%
%%Fix some choice of $j_2,\dots,j_n$, and let $G_i := G_{i,j_i}$.
%
%Then, by \cite[p. 136]{fis01}, for each $i$ there is a
%positive natural number
%$k_i$ and a function $\phi_i$ holomorphic on a disk $D_i$ centered
%at $0$ so that, for some open set $U_i$ in $\C^2$ containing the
%point
%%$(P_1,P_i)$, 
%$$
%\{ (z_1,z_i)\, \in\, U_i
% \ : \ G_i(z_1,z_i) = 0 \} \= \{ ( P_1 + t^{k_i}, \phi_i(t) )
%\ : \ t \, \in \, D_i \} 
%$$
%and so that the map
%%$$
%t \ \mapsto \ (t^{k_i}, \phi_i(t) )
%$$
%is injective on $D_i$.
%
%Let $k = {\rm lcm}(k_i)$. 
%Define $\Phi$ on $D$, the intersection of the disks $D_i$, by
%\be
%\label{eq57}
%\Phi \ : \ t \ \mapsto \ (t^k, \phi_2 (t^{k/k_2}), \dots, \phi_n(
%t^{k/k_n}) )  .
%\ee
%Then $\Phi$ is one-to-one on $D$ and maps onto a neighborhood of
%$P$ in the irreducible component chosen.
%
The disk $\Delta$ and the function $h$ depend on the choice
of point $P$ and the choice of irreducible component.
For each such choice, take the disjoint union of $S'$ with $\Delta$, and
then identify the points of $\Delta\setminus \{ 0 \}$ with their image
in $S$ under $h$. This new Riemann surface is the desired $S$.
\vs
(iii)
Define $\psi = \tilde{h}^{-1} \circ h$. This is a well-defined
biholomorphic bijection from $S' = S\setminus E$ onto $\tilde S'
= \tilde S\setminus \tilde E$.

Now consider $\psi$ on a punctured neighborhood of a point $Q$ in $E$.
In local coordinates, it is a bounded analytic function on a
punctured disk, therefore it extends to be holomorphic at $Q$.
Moreover, by continuity, as $\zeta $ tends to $Q$,
$\tilde{h}^{-1} \circ h (\zeta)$ will converge to some point in $\tilde
S$, since $\tilde S$ is unpunctured. 

Extend $\psi$ in this manner to all of $E$. The new $\psi$ will
remain injective, since otherwise there would be points $\tilde Q$
in $\tilde S$  for which every neigborhood contained 2 (or more)
disks intersecting only at $\tilde Q$. 
Moreover, by applying the same argument to $\psi^{-1}$, we see that
$\psi$ must also be surjective.
\ep

%\rem
%The only place in the paper where we use explicitly the properness of holomaps is
%in the assertion in the proof of (i) that $h(S)$ is necessarily relatively closed in $U$.
%If we are willing to add this as a hypothesis in Theorem~\ref{thmb1} (i),
%then we can remove ''proper" from  the definition of holomap.

\section{Cofinite algebras}~\label{secc}

Cofinite sub-algebras of uniform algebras were studied in depth by T. Gamelin, 
who showed that they can be obtained inductively by imposing one linear restriction
at a time on a cofinite algebra of lower codimension \cite{gam68}.
He stated his results in the context of uniform algebras, but his proofs of the
parts that we shall use (Theorem~\ref{thm2} below), as he
remarks in the paper, carry over to our setting.

%Let $B$ be a uniform algebra, and $M_B$ denote its maximal ideal space.
%A point derivation of $B$ at $\theta \inn M_B$ is a linear functional $D$ on $B$ that
%satisfies
%$$
%D(fg) \= \theta (f) Dg + \theta(g) Df .
%$$
%A subalgebra $H$ of $B$ is called a $\theta$-subalgebra
%if there is a finite sequence of algebras
%$B = B_0 \supersetneq B_1 \supersetneq \dots \supersetneq B_k = H$ such that
%each $B_j$ is the kernel of a continuous point derivation of $B_{j-1}$ at $\theta$. 

Before giving his result, let us make two more definitions.
\bd
\label{def21}
A linear functional on $O(S)$ is called local if it comes from a
finitely supported distribution, \ie is of the form
$$
\Lambda(f) \= \sum_{i=1}^m \sum_{j=0}^{n_i} a_{ij}
f^{(j)}(\alpha_i) .
$$
\ed

\bd
\label{def22}
A connection on $\{ \a_1,\dots , \a_m \} \subset S$
is a finite dimensional set $\Gamma$ of local functionals $\Lambda$ supported in
$\{ \a_1,\dots , \a_m \}$. We say $\Gamma$ is algebraic if
$\Gamma^\perp := \{ f \ \in \ O(S)\ : \ \Lambda(f) = 0 \ \forall \,
\Lambda\, \in \, \Gamma \}$ is an algebra.
We say $\Gamma$ is irreducible if every function $f$ in
$\Gamma^\perp$ satisfies $f(\a_1) = f(\a_2) = \dots = f(\a_m)$.
\ed
\vs
Note that if the support of $\Gamma$ is a singleton, it is
by definition irreducible.
In \cite{bis58}, E.~Bishop conidered analogous objects to our Definitions~\ref{def21}
and \ref{def22}. In his terminiology, each $\Lambda$ in the algebraic
connection $\Gamma$
is a {\em homogeneous linear operator of order no more than $\max (n_i)$ for the algebra
$\Gamma^\perp$}.

\bt
\label{thm2} [Gamelin]
(i) If $\Gamma$ is an algebraic connection on a Riemann surface $S$
then $\Gamma^\perp$ is a cofinite algebra.

(ii) Every cofinite subalgebra  of $O(S)$ is of the form
$\Gamma^\perp$ for some algebraic connection $\Gamma$.

(iii) Every algebraic connection is the union of a finite number of
irreducible algebraic connections with disjoint supports.

(iv) Every cofinite algebra $A$ of $O(S)$ can be obtained by a chain
$ A \= A_n \subsetneq A_{n-1}  \dots  \subsetneq A_1 = O(S)$, where 
each $A_{j+1}$ is an algebra of codimension one in $A_j$ and is obtained 
either as the kernel of a point derivation on $A_j$ or by identifying two points of the maximal
ideal space of $A_j$.
\et
\bp
Part (i) is obvious. Parts (ii) - (iv) were proved by T.~Gamelin 
\cite[Lemma 9.7 and Thm. 9.8]{gam68}.
\ep
\vs
Before giving our holization theorem, we need some lemmata.
\bl
\label{lem1}
Let $N, K \geq 1$ be coprime.
Let $h_1(z) \=  z^N + \sum_{i= N+1}^\i a_i z^i$ and
$h_2(z) \=  z^K + \sum_{i= K+1}^\i a_i z^i$ be holomorphic on some
disk centered at the origin of radius $ R > 0$.
Then there is some $r >0$ such that, in the disk algebra $A(\D_r)$,
the algebra generated by $h_1$ and $h_2$ contains all
holomorphic functions that vanish at the origin to a high order
(that depends only on $N$ and $K$).
Moreover, for all large $L$, there is some function $F_L$
holomorphic on $\D_r \times \D_r$, such that $F_L (h_1(z),h_2(z) )
= z^L$.
\el
\bp
Observe that the leading terms of the functions $h_2^j$, the
monomials $z^{jK}$, have powers that are all different mod $N$ for
$1 \leq j \leq N$. Let $\sigma$ be a permutation of
 $\{0,1,\dots,N-1\}$ such that
the leading term of $h_2^{\sigma(j)}$ has a power congruent to $j$
mod $N$.

Define the functions $f_j$ by
$$
f_j (z) \= (h_1)^{n_j} h_2^{\sigma(j)} ,
$$
where the non-negative integers $n_j$ are chosen so that 
the powers of the leading terms of the functions $f_j$ form the
arithmetic progression $\{ mN,mN+1,\dots, (m+1)N -1\}$.

We wish to obtain the functions $z^{mN}$ and $z^{mN +1}$ in the
algebra generated by $h_1$ and $h_2$.
Once we do this, the algebra generated by $z^{mN}$ and $z^{mN +1}$
contains each of the functions
$$
(z^{mN})^{mN-j} (z^{mN +1})^{j} \= z^{m^2N^2+j}
$$
for $0 \leq j \leq mN$, and so contains every monomial $z^k$ for
$k \geq m^2 N^2$.

Starting with $f_0$, one can subtract a multiple of $f_1$ 
to cancel the coefficient of $z^{mN+1}$, and then multiples of
$f_2, \dots, f_{N-1}$ to cancel the next few coefficients. Then one
subtracts multiples of $h_1^{i} f_j$ to cancel out the coefficients
one by one. If we can show that this process converges, we are
done.

Without loss of generality, $R$ can be taken greater than $1$,
and so
there is some constant $C$ that bounds the
Taylor coefficients of each $f_j$.
Let us relabel the functions $$
\{ f_1, f_2, \dots, f_{N-1}, h_1 f_0, \dots , h_1 f_{N-1}, h_1^{2} f_0,
\dots \}
$$
as $\{ g_1, g_2, \dots \}$. Then at the $k^{\rm th}$ stage of the
cancellation, we have the function
$$
f_0\ -\ \sum_{i=1}^k a_i g_i ,
$$
and at the $(k+1)^{\rm st}$ stage we have the function
$$
f_0\ -\ \sum_{i=1}^k a_i g_i \ - \ \left(\hat{f_0}(mN +k+1) -
\sum_{i=1}^k a_i \hat{g_i} (mN+k+1) \right) g_{k+1}.
$$
So we have the inequality
\be
\label{eq34}
|a_{k+1} | \ \leq \ C ( 1 + \sum_{i=1}^k |a_k | ).
\ee
It follows from (\ref{eq34}) by induction that
$$
|a_k | \ \leq \ C\,(C+1)^{k-1} \qquad \forall\ k .
$$

On the disk of radius $r$,
the modulus of $g_k$ is bounded by a constant times $r^{k/N}$.
Therefore if $r$ is chosen so that
$$
r \ < \ (C+1)^{-N} ,
$$
then the series $\sum a_k g_k$ will converge absolutely, and so the
function
$z^{mN}$
will be in 
the closure in $A(\D_r)$ of the algebra generated by $h_1$ and $h_2$,
as desired.
Similarly, starting with $g_1$, one gets the function $z^{mN+1}$.

Finally, the function $F_L$ is a finite product of functions of 
the form $\sum b_{ij} z_1^i z_2^j$, where 
$$
|b_{ij}|\  \leq \ (C+1)^{iN}
$$
and, for each $i$, there are at most $N$ distinct $j$'s with
$b_{ij} \neq 0$. Therefore $F_L$ is analytic on $\D_r \times \D_r$,
as desired.
\ep
\bl
\label{lem2}
Let $T$ be an additive subsemigroup of $\NN$, and suppose that
${\rm gcd} \{ t\ \in \ T \} = k$. Then $k \NN \setminus T$ is finite.
\el
\bp
By hypothesis, there exist $t_1,\dots,t_n$ in $T$ such that
${\rm gcd}(t_1,\dots,t_n) = k$. Let $\tau = {\rm lcm}(t_1,\dots,t_n)$.
Then, since the subgroup generated by $T$ is all of $k \Z$, we have
integers $m_{ij}, 1 \leq i \leq \tau/k,\ 1\leq j \leq n$, so that
$$
\sum_{j=1}^n m_{ij} t_j \= ik .
$$
Let $L \in \NN$ be such that $L+m_{ij} \geq 0$ for every $i,j$.
Then $T$ contains the set
$$
\{ \sum_{j=1}^n t_j (m_{ij} + (L+r) \frac{\tau}{t_j} )\ :
\ 1 \leq i \leq \tau/k,\
r \in \NN \} ,
$$
and this set equals $$
\{ s \in \, k\NN \ : \ s > L \tau \} .
$$
\ep
\bd
Let $f$ be a holomorphic function in a neighborhood of $0$.
By the {\em order} of $f$, written ${\rm ord}(f)$,
we mean the smallest number such that the corresponding Taylor
coefficient for $f$ is non-zero.
\ed

\bl
\label{lem3}
Let $h = (h_1,\dots, h_n)$ be an $n$-tuple of analytic functions
on $\D$, with $h_i(0) = 0$ for each $i$.
Suppose $h_1(z) = z^N$, and that the $n$-tuple separates points 
on $\D$. Let $k > 1$ be a divisor of $N$. Then the algebra
generated by $h_1,\dots,h_n$ contains a function 
$f$ with ${\rm ord} (f)$ not a multiple of $k$.
\el
\bp
Suppose not. Then every function in $\A$, the algebra
generated by $h_1,\dots,h_n$, has order a multiple of $k$. Some
function must have some Taylor coefficient that is not a multiple of
$k$, as the algebra separates points.
Choose $f$ in $\A$ such that the difference between the first power
of $z$ that is not a multiple of $k$ and whose Taylor coefficient
is non-zero and the order of $f$ is minimal. So 
$$
f(z) \= \sum_{i = i_0}^\infty a_{ik} z^{ik} + a_r z^r +
O(z^{r+1}) ,
$$
where $r$ is not a multiple of $k$ and $r - i_0 k$ is minimal.
Let $N = nk$.
Now consider
$$
g(z) \= f^n(z) - z^{i_0 N} .
$$
Then ${\rm ord}(g) \geq i_0 N + 1$, and it contains a term whose
$(n-1)i_0 k +r$ Taylor coefficient is non-zero.
The difference between this and the order of $g$ is less than or
equal to $ r - i_0 k -1$, a contradiction.
\ep

It is a theorem of H.~Rossi \cite{ros61} that if $X$ is a compact set in $\C^n$,
then $R(X)$ can be generated by $n+1$ functions. So the condition in (iii) below is not too restrictive. 
It clearly applies, for example, to any smoothly bounded planar domain.

\bt
\label{thmb3}
(i) If $h : S \to V \subseteq U$ is a holomap from a Riemann
surface onto  a \gas $V$, then 
$$
A_h \ := \ \{ F \circ h \ : \ F \ \in \ O(V) \}
$$
is a cofinite subalgebra of $O(S)$.

(ii) Conversely, if $A$ is a cofinite subalgebra of $O(S)$ for some
good Riemann surface $S$, then there
exists a holomap $h : S \to \C^n$ such that $A = A_h$.

(iii) Let  $V$ be a \hacp, and assume that $O(V)$ can be generated by a 
finite number of rational functions that are continuous on $\overline{V}$.   
If $A$ is a cofinite subalgebra of $O(V)$,
then there is a \hac $W$ and a holomap $h$ from $V$ onto $W$ such that
$A = A_h$.
\et
%We say such an $h$ is a {\em holization} of $A$.

\bp
(i)
(a)
Let $g \in O(S)$. Then $g \in A_h$ if $F = g \circ h^{-1}$ is in
$O(V)$. By Theorem~\ref{thmb1} (iii), we can assume $S$ is the
surface constructed in the proof of Theorem~\ref{thmb1} (ii).
Then $F$ is analytic on $V \setminus B$. 

For $F$ to be analytic at the point $P$ in $B$, we need certain
consistency conditions on $g$ and its derivatives on the pre-images
of $P$. Let $Q$ be a pre-image of $P$. We show that if $g$ vanishes
to high enough order at $Q$, then there will be an analytic
function $F$ on a neighborhood of $P$
such that $F \circ h$ agrees with $g$ on a
neighborhood of $Q$.

Choose a disk $\Delta$ containing $Q$ such that 
$\overline{\Delta} \cap E = \{ Q \}$, and let $z$ be a local
coordinate at $Q$ so that $\Delta = \{ z : |z| < 1 \}$,
and such that $h_1(z) = z^N$ (after, perhaps, permuting the
$h_i$'s).

Consider the subalgebra of $O(\Delta)$ generated by $h_1, \dots,
h_n$. Let $T$ be the set of natural numbers $t$ such that there is
a function in the algebra generated by $h_1, \dots,
h_n$ whose first non-zero Taylor coefficient is the $t^{th}$.
Then $T$ is a semigroup, and 
so by Lemma~\ref{lem2}, is $k \NN$ minus, perhaps, a finite set.

By Lemma~\ref{lem3}, we must have $k=1$. Therefore $T$ contains
numbers that are coprime, so by Lemma~\ref{lem1}, 
$\{ F \circ h \ : \ F \ \in \ O(P) \}$ will contain all 
functions analytic on a neigborhood of $Q$ that vanish at $Q$ to
high order.

\vs
(b)
In part (a), we showed that for each pre-image $Q$ of $P$, any
function $g$ that vanishes at $Q$ to high order
can be locally represented as $F_Q \circ h$ for some 
$F_Q$ analytic near $P$. To finish the proof, we need to show that
we can isolate the (finitely many)
different sheets at $P$. This means that if
$Q_1$ and $Q_2$ are distinct preimages of $P$, and, as in the proof
of part (a) we have disks $\Delta_1$ and $\Delta_2$ with maps
$$
h^i \ : \ z_i \, \mapsto \, (z_i^{N_i}, h_2^i(z_i), \dots ) \qquad
i=1,2,
$$
then there is some function $F$ analytic near $P$ such that
$F \circ h^1(z_1) $ is a power of $z_1$ and
$F \circ h^2$ is identically $0$. Indeed, by part (a) we could then
find, for every set of functions $g_1, \dots, g_r$ analytic on
neighborhoods of $Q_1, \dots, Q_r$ in $S$ and vanishing to high
enough order, functions $F_i$ such that the function
$F_i \circ h$ will be $g_i$ near $Q_i$ and $0$ near all other
$Q_j$'s. The sum of the $F_i$'s would then be the desired function.

By the proof of \cite[Thm. III.B.19]{gur} (applied to the map that
is projection
onto the first coordinate),
$h^2(\Delta_2)$ is a one dimensional analytic subvariety, and there
exists some open set $U'$ containing $P$ with the property that
for every 
point $R$ in $U' \setminus h^2(\Delta_2)$ there is some function
$G$ in $O(U')$ with the property that $G(R) \neq 0$ and
$G |_{h^2(\Delta_2)} \equiv 0 $.

Choose $R$ in $h^1(\Delta_1)\setminus \{ P \}$.
Then we get a $G$ such that $G \circ h^1$ is non-constant, and 
$G \circ h^2 \equiv 0$. By shrinking the neighborhood $\Delta_1$
further, if necessary, we can assume that the only zero of $G \circ
h^1$ on $\overline{\Delta_1}$ is at $0$. 
Therefore $G \circ h^1 (z_1)$ is  $ z_1^N$ times some unit in 
$O(\Delta_1)$. By part 
(a), $\{ F \circ h^1 \ : \ F \, \in \, O(P) \}$ contains $z_1^k 
O(\Delta_1)$ for some $K$, so the set 
$\{ (F\cdot G) \circ h \ : \ F \, \in \, O(P) \}$
will contain all functions 
on $\Delta_1$ that vanish to order $N+K$ at $Q_1$, and will be
identically zero on  a neighborhood of $Q_2$.

\vs
(ii) It follows from Theorem~\ref{thm2} that there is a natural
number $N$ with the property that every function in $O(S)$ that
vanishes at every point of the connection to order $N$ is in $A$.
%So it is sufficient to prove that this algebra is finitely
%generated as a Fr\'echet algebra.

Since $S$ is good, there is a non-singular holomap
$f$ from $S$ into 
some bounded set $U \subseteq \C^k$ such that
polynomials in $f_1, \dots, f_k$ are dense in $O(S)$.
Let $V = f(S)$, and let the image of
the points of the connection under $f$ be
be $\{ \alpha_1, \dots, \alpha_m \}$; without loss
of generality assume no $\alpha_i$ is $0$. 
Then the polynomials in
$\C[z_1,\dots, z_k]$ that vanish to order $N$ 
at each $\alpha_i$ are finitely generated, as an algebra, by
polynomials of the form
\be
\label{eqd8}
(z_1 - \alpha_1^1)^{j_{11}} (z_2 - \alpha_1^2)^{j_{12}} \cdots
(z_k - \alpha_m^k)^{j_{mk}} \, z_1^{i_1} \cdots z_k^{i_k},
\ee
where $j_{11} + j_{12} + \cdots + j_{mk} = N$ and $0 \leq i_1,
\dots, i_k \leq N-1$.

If necessary, add to the collection of 
functions given in (\ref{eqd8}) a finite
number of additional polynomials that vanish to order less than $N$
and satisfy the conditions of the connection.
Then we will have a finite number of polynomials
$g_1,\dots, g_n$ with the property that the algebra generated
by these
functions will give all polynomials that satisfy the connection on
$V$.
Therefore the functions $h_1 := g_1 \circ f, \dots, h_n := g_n \circ f$ are
functions in $A$ with the property that the algebra they generate
is the intersection of $A$ with the algebra of polynomials in $f$.
Since polynomials in $f$ are dense in $O(S)$ by hypothesis,
any function in $A$ is the u.c.c. limit of a sequence of
polynomials in $h = (h_1, \cdots , h_n)$. Since $V$ is a non-singular
analytic subvariety of $U$, 
we have therefore that $A = A_h$.

Finally we must prove that $h$ is proper, and therefore a holomap.
But this follows because $h$ is one-to-one off the set $\{ \alpha_1,\dots,\alpha_m \}$,
and $h$ maps open subsets of $S$ to relatively open subsets of $h(S)$.

\vs
(iii) Let $V$ be a \hac as in the statement of the theorem, with $V$ contained in the algebraic
curve $\CC$. Assume that $\CC$ has no irreducible components disjoint from $V$.
We can choose $f: V \to U \subseteq \C^k$ as in part (ii), with the additional property 
that the $f_i$'s
are rational and bounded on $V$.
Moreover, $f$ extends to be a regular rational function on $\CC$ less, at most, a finite set.
Choose the $g_j$'s as in part (ii), and let $h = g \circ f$.
By Chevalley's constructibility theorem \cite[Thm. 3.16]{har92}, 
$h(\CC)$ is contained in some algebraic curve ${\cal D} \subseteq \C^n$.
As $h$ is a homeomorphism off the connection, $ W = h(V)$ is relatively open in
${\cal D}$, so $W$ is a \hacp, as required.
\ep
%Finally we must show that $f$ is a holomap.
%Since $A$ is a cofinite algebra, $f$ can have only a finite number
%of singular points. To see that $f : S \to g(U)$ is proper,
%Note that $g : U \to \C^n$
%is a polynomial map, 
%\rem
%The fact that we are working with the whole algebras $O(S)$ means that the boundary behavior
%can largely be ignored. When using the normal Banach spaces of functions on $S$,  
%it is a much more delicate issue to decide when a pair of functions generates a finite codimensional
%subspace. See \eg \cite{stth03} and \cite{wer58}. 

\section{Isomorphism theorem}\label{seccc}

We now reach the pay-off for setting things up carefully, the isomorphism theorem,
which allows us to pass between algebraic and geometric information.

If $\Lambda$ is a linear functional on $O(S_1)$, and $\psi: S_1 \to S_2$ is 
a holomoprphic map, we define the push-forward $\psi_* (\Lambda)$ by
$$\psi_* (\Lambda) [f] \=  \Lambda[ f \circ \psi] .$$
Applying this to each element of a connection $\Gamma$ on $S_1$, one gets the
push-forward $\psi_* (\Gamma)$ on $S_2$.

\bt
\label{thmcc1}
Let $S_1$ and $S_2$ be Riemann surfaces, and, for $r=1,2$, let $h_r : S_r \to V_r \subseteq U_r$ be a holomap
onto a \gasp. Then the following are equivalent:

\noindent
(i) There is a biholomorphic homeomorphism $\phi : V_1 \to V_2$.

\noindent
(ii) There are biholomorphic homeomorphisms $\psi : S_1 \to S_2$ 
and $\phi : V_1 \to V_2$ such that
$h_2 \circ \psi = \phi \circ h_1$.

\noindent
(iii) There is a biholomorphic homeomorphism $\psi : S_1 \to S_2$ such that
\beq
A_{h_1} (S_1) &\ :=\ & \{ F \circ h_1 \ : \ F \inn O(V_1) \} \\
&=& %A_{h_2 \circ \psi} \ :i
\{ G \circ h_2 \circ \psi \ : \ G \inn O(V_2) \} .
\eeq

\noindent
(iv) There is a biholomorphic homeomorphism $\psi : S_1 \to S_2$ that pushes forward the
connection on $S_1$ induced by $h_1$ to the connection on $S_2$ induced by $h_2$.
\et
\bp
(i) $\Rightarrow$ (ii): Let $E_1$ and $E_2$ be the finite subsets of $S_1$ and $S_2$ respectively,
off of which $h_1$ and $h_2$ are injective and non-singular.
Let $B_1 = h(E_1)$ and $B_2 = h(E_2)$.
Define $\psi$ by
$$
\psi \= h_2^{-1} \circ \phi \circ h_1 .
$$
Then $\psi$ is holomorphic on $S_1 \setminus (\phi \circ h_1)^{-1} (B_2)$.
For any $Q$ in $ (\phi \circ h_1)^{-1} (B_2)$, the function $\psi$ maps a small
punctured neighborhood of $Q$ to a small punctured neighborhood in $S_2$ that is
disjoint from $E_2$. So the singularity at $Q$ is removable, and $\psi$ extends to 
be holomorphic from $S_1$ to $S_2$. 
The same argument shows that $\psi^{-1}$ is holomorphic.
%By construction, $\psi$ is proper and one-to-one off a finite set.

(ii) $\Rightarrow$ (iii): Define $G = F \circ \phi^{-1}$.

(iii) $\Rightarrow$ (iv): Let the connections $\Gamma_1,\Gamma_2$
be given by
 $$\Gamma_r = A_{h_r} (S_r)^\perp,\qquad r =1,2.
$$
Define $\tau_\psi : O(S_2) \to O(S_1)$ by $\tau_\psi (g) = g \circ \psi$.
Then (iii) says
\beq
\Gamma_1 
&\=& \{ \Lambda \inn O(S_1)^\ast \ : \ \Lambda ( F \circ h_1) = 0 \ \forall \
F \inn O(V_1) \} \\
&\=& \{ \Lambda \inn O(S_1)^\ast \ : \ \Lambda \circ \tau_\psi ( G \circ h_2 ) = 0 \ \forall \
G \inn O(V_2) \} \\
&\=& \Gamma_2 \circ \tau_\psi^{-1} .
\eeq
Therefore 
$$
\Gamma_1 \circ \tau_\psi \= \Gamma_2 .
$$

(iv) $\Rightarrow$ (i): 
%Define $\phi = h_2 \circ \psi \circ h_1^{-1}$.
%Then $\phi$ is holomorphic off $B_1$. 
%Let $P$ be a point in $B_1$. Locally, each sheet at $P$ is mapped to a 
%sheet in $V_2$, and the sheets meet at a unique point $Q$, so $\phi$
%extends to be continuous at $P$.
% and holomorphic on each sheet separately.
For every $G$ in $O(V_2)$, the function $G \circ h_2 \circ \psi$
is in $O(S_1)$. Moreover, if $\Lambda$ is in $\Gamma_1$, then
$$
\Lambda(G \circ h_2 \circ \psi) \=
(\Lambda \circ \tau_\psi ) (G \circ h_2) \= 0 .
$$
Therefore $G \circ h_2 \circ \psi$ is in $A_{h_1}$ for every $G$.
By Theorem~\ref{thmb3}, this means there always exists some $F$ in $O(V_1)$
so that
$G \circ h_2 \circ \psi = F \circ h_1$.

Letting $G$ run over the coordinate functions, we get that there is a holomorphic mapping $\phi$
on $V_1$ such that $$
h_2 \circ \psi \= \phi \circ h_1 . $$ 
As $\phi$ is one-to-one by construction, the same argument applied to $\psi^{-1}$ 
yields that $\phi$ is a biholomorphic equivalence.
\ep

%The proof follows from Theorems~\ref{thmb1} and ~\ref{thmb3}.

We shall say that two \gass $V_1$ and $V_2$ are isomorphic 
if there is a biholomorphic homeomorphism between them. Note that even if they are both
\hacsp, one cannot assume that the map is algebraic; for example, a disk and a square
in the plane are isomorphic by our definition.
We shall say that two cofinite algebras are isomorphic if they are related as in (iii)
of Theorem~\ref{thmcc1}.

A very similar argument to the proof of Theorem~\ref{thmcc1}
allows us to tell when there is any holomorphic map between two 
\gassp. Using the same notation, we have:
\bt
\label{thmcc2}
Let $h_r : S_r \to V_r$ be  holomaps
onto  \gassp. Then the following are equivalent:

\noindent
(i) There is a holomorphic map $\phi : V_1 \to V_2$.

\noindent
(ii) There are holomorphic maps $\psi : S_1 \to S_2$ 
and $\phi : V_1 \to V_2$ such that
$h_2 \circ \psi = \phi \circ h_1$.

\noindent
(iii) There is a holomorphic map $\psi : S_1 \to S_2$ such that
$$
A_{h_1} (S_1) \ \supseteq\  \{ A_{h_2}(S_2) \, \circ \psi \} .
$$

\noindent
(iv) There is a holomorphic map $\psi : S_1 \to S_2$ such that 
$$
\Gamma_1 \circ \tau_\psi \ \subseteq \ \Gamma_2.
$$
\et

\section{Generators}
\label{secy}

If $A$ is a cofinite subalgebra of $O(S)$, and $f_1,\dots,f_n$ are functions in $A$,
it is of interest to know when they generate the whole algebra.
Let $\alg(f_1,\dots,f_n)$
denote the closure of the polynomials in $f_1,\dots, f_n$ (in the Fr\'echet topology).

In order for a unital sub-algebra $A$ to equal $O(S)$, the following two conditions are clearly
necessary:

(SP): The functions in $A$ separate points of $S$. 

(ND): At every point of $S$ there is a function in $A$ whose derivative does not vanish there.

We shall show that these conditions are sufficient in two particular cases.
The first condition, named after J.~Wermer, is topological; the second, named after E.~Bishop,
is a Mergelyan-like theorem, in which the algebra $A$ must take the topology of $S$ into account.

(W): Suppose that $S$ is the interior of a compact set $D$ in a Riemann surface ${\cal S}$,
and the boundary of $S$ is a simple closed analytic curve.

(B): Suppose that there is a sequence of compact sets $K_n$ such that 
each $K_n$ is contained in the interior of $K_{n+1}$, the union of all the $K_n$ is
$S$, and there exists $n_0$ such that, for all $n \geq n_0$, for all $p$ in $S \setminus K_n$,
there exists some $f$ in $A$ with $|f(p)| > \| f \|_{K_n} $.

\bt
\label{thmy1}
Let $A$ be a unital sub-algebra of $O(S)$ satisfying (SP) and (ND). If either (W) or
(B) is also satisfied, then $A$ is dense in $O(S)$.
\et
\bp
First assume (W). Then one can find an increasing sequence of compact sets $D_n$ whose union is
$S$ and such that each $D_n$ has a boundary that is a simple closed analytic curve $\gamma_n$.
Moreover, one can find a function $\phi$ in $A$ whose derivative is non-zero anywhere
on $\gamma_n$. (Indeed, choose $\phi$ to be any function whose derivative does not vanish identically.
Then its zeroes are isolated, and one can perturb $\gamma_n$ to miss them all).

By a result of Wermer \cite[Appendix]{wer58}, the algebra $A$ is then uniformly dense
in $A(D_n)$, the functions analytic on the interior of $D_n$ and continuous up to the boundary.
Now suppose that $A$ were not dense in $O(S)$. Then by the Hahn-Banach theorem, there would be a 
non-zero linear functional $\Lambda$ on $O(S)$ that annihilated $A$. Every element
of the dual of $O(S)$ comes from integration against a compactly supported measure $\mu$
\cite[Thm. 4.10.1]{ed95}.
For $n$ large enough, $\mu$ will be supported on $D_n$. But as integration against $\mu$
annihilates a dense subspace of $A(D_n)$, it must annihilate all holomorphic functions on $S$
restricted to $D_n$, a contradiction.

Now assume (B). Again, argue by contradiction, and assume that there is a compactly supported measure
$\mu$ that annihilates $A$ but not all of $O(S)$. Let the support of $\mu$ be $K$. Let 
$\hat K$ be $\{ p \inn S \ : \ \forall \ f \inn A, \  |f(p)|  \leq \| f \|_K \}$.
By assumption (B), $\hat K$ is compact. By a result of Bishop \cite{bis58}, $A$ is dense
in $A(\hat K)$, and again we achieve a contradiction.
\ep

It is not true that every subalgebra of $O(S)$ is contained in a proper 
cofinite subalgebra: consider, for example, the algebra of polynomials
inside the algebra of all holomorphic functions on an annulus. Theorem~\ref{thmy1} shows, however,
that if the set $S$ satisfies (W), or the algebra $A$ satisifies (B), then it is either dense, or contained
in a codimension one subalgebra of $O(S)$.

Given a particular cofinite subalgebra, we know from Theorem~\ref{thmb3}
that it can be holized into $\C^n$ for some $n$. To prove that $n$ can be chosen 
to be small, we need to show that some choice of $f_1,\dots,f_n$ in $A$ do not generate a proper subalgebra.
The next theorem says that, in certain important circumstances
(such as trying to holize an algebra on the disk into $\C^2$) you only need to worry about the obvious obstructions.

\bt
\label{thmy2}
Let $S$ be a set satisfying (W).
Let $\Gamma$ be a connection on $S$, and let $A = \Gamma^\perp$. Let $F= \{\al_1,\dots,\al_m\}$
be the support of $\Gamma$. Let $f = (f_1,f_2)$ be a pair of functions in $A$.
Assume that

\noindent
(i) If $f$ identifies a pair of points in $S$, then both points are in $F$.

\noindent
(ii) The derivative of $f$ is non-vanishing on $S \setminus F$.

Then
$\alg(f_1,f_2)$ is a cofinite algebra given by a connection that is supported on $F$.
\et
\bp
Choose an increasing sequence $D_n$ as in the proof of Theorem~\ref{thmy1}, and moreover
assume that $F \subset {\rm int}(D_1)$. If necessary, perturb the boundaries of the $D_n$'s
so that either $df_1$ or $df_2$ is never zero on $\partial D_n$.
Then by the proof of Theorem~1.2 in \cite{wer58}, there is an integer $k$ that is independent of
$n$ and such that every function in $A(D_n)$ that vanishes to order $k$ on $F$ is uniformly approximable
by functions in $\alg(f_1,f_2)$. Therefore $\alg(f_1,f_2)$ contains all functions on $S$ that
vanish to order $k$ on $F$, so by Theorem~\ref{thm2} we are done.
\ep

Because of the topology we are using, many subtleties of boundary behavior disappear.
In \cite{stth03}, Stessin and Thomas consider the problem of determining the $H^p$-closure of the algebra generated
by a pair of functions.

If an algebra is generated by an $n$-tuple $h = (h_1,\dots,h_n)$, then $h$ is not necessarily a
holomap: if $h(\partial S)$ is not disjoint from $h(S)$, then one cannot find an open set $U$
containing $h(S)$ such that $h:S \to U$ is proper.
Conversely, if $A = A_h$ for some holomap $h: S \to V \subset U$, then the polynomials in $h$ need not be
all the functions in $A$, if there are holomorphic functions on $V$ that are not approximable
by polynomials. If all the functions $h_i$ are inner functions that are continuous up to the boundary
of $S$, however, then everything is nice.
Indeed, the map $h$ is then a a holomap into the polydisk $\D^n$.
By Theorem~\ref{thmb1}, the set $V = h(S)$  is a \gasp.
By 
Cartan's theorem~\cite{car51}, any holomorphic function $F$ on $V$ extends to
a holomorphic function on $\D^n$, and is therefore uniformly approximable on compact sets
by polynomials. Thus we have:
\bprop
\label{propy1}
Let $S$ be a smoothly bounded Riemann surface. Let $h = (h_1,\dots, h_n)$ be an $n$-tuple of
holomorphic functions that are continuous up to $\partial S$ and have modulus one everywhere on
$\partial S$. Then 
$\alg(h_1,\dots,h_n) \= A_h$.
\eprop

\section{Pairs of Blaschke products}\label{secbp}

It is natural to ask if a given cofinite subalgebra of $O(\D)$ can be holized by a pair of Blaschke
products.
S.I.~Fedorov proved that on every $C^\omega$-bounded planar domain $\Omega$ there 
are a pair of continuous inner functions $\phi_1$ and $\phi_2$ that continue analytically across the boundary
and such that the pair $\Phi = (\phi_1,\phi_2)$
separate points of $\Omega$ and $d\Phi \neq 0$ on $\Omega$ \cite{fed91}. This means that $\Phi$ holizes
$\Omega$ as a distinguished variety. (The fact that $\phi_1$ and $\phi_2$ satisfy an algebraic equation
was proved by Fedorov by extending them to the Schottky double, but also follows from general arguments
\cite{agmc_dv}.) However, W.~Rudin showed that if the connectivity of $\Omega$ is greater than
one, then $\Phi$ cannot be one-to-one on $\partial \Omega$ \cite{rud69b}.

We shall want to know how many multiple points a pair $(f,g)$ of finite Blaschke products (each
of degree
greater than one)
has on $\C^2$, \ie how many solutions $(\l,\mu) \inn \C^2$ there are to the simultaneous equations
\se\att
\begin{eqnarray}
\nonumber
f(\l) &\=& f(\mu) \\
g(\l) &\=& g(\mu) .
\label{eqbp1}
\end{eqnarray}
When counting these, we shall consider $(\l,\mu)$ and $(\mu,\l)$ to be distinct solutions, and we shall
need to count solutions with multiplicity.
To this end, write 
\se\att\begin{eqnarray}
\label{eqbp2}
f(z) &\=& \frac{p^\sim(z)}{p(z)} \\
\att
\label{eqbp3}
g(z) &\=& \frac{q^\sim(z)}{q(z)} 
\end{eqnarray}
where $p(z) = \prod_{i=1}^m (1- \overline{\al_i}z),\ p^\sim(z) =  \prod_{i=1}^m (z - \al_i)$
and
$q(z) = \prod_{i=1}^n (1- \overline{\beta_i}z),\ q^\sim(z) =  \prod_{i=1}^n (z - \beta_i)$.
%Define the multiplicity of the solution to (\ref{eqbp1}) at the point $(\l,\mu)$ as the 
%intersection multiplicity  of the curves $\CC_F$ and $\CC_G$
%at the point  $(\l,\mu)$, where  
Let $\CC_F$ and $\CC_G$ be the zero sets of the polynomials
\se\att
\begin{eqnarray}
%\CC_F &\=& \{ 
F(z,w) &\ := \ & \frac{p(z)p^\sim(w) - p(w)p^\sim(z)}{z-w}
\label{eqbp4} \\
\att
%\CC_G &\=& \{ 
G(z,w) &\ := \ & \frac{q(z)q^\sim(w) - q(w)q^\sim(z)}{z-w}
\label{eqbp5} 
\end{eqnarray}
respectively.

Note that generically $F$ in (\ref{eqbp4}) is a polynomial of degree $2(m-1)$; the
coefficient of $z^{m-1}w^{m-1}$ is 
\be
\label{eqbp55}
(-1)^{m+1}  \left[ \sum_{i=1}^m \al_i \prod_{j=1}^m \bar \alpha_j - \sum_{i=1}^m \prod_{j \neq i} \bar \alpha_j \right] 
\=
(-1)^{m+1} \overline{f'(0)} .
\ee
Let $\hat F(t:z:w)$ and $\hat G(t:z:w)$
be the homogeneous polynomials of degrees exactly 
$2(m-1)$ and $2(n-1)$ so that
\beq
\hat F(1:z:w) &\=& F(z,w) \\
\hat G(1:z:w) &\=& G(z,w).
\eeq

%If $\{ \al_i \}$ is such that (\ref{eqbp55}) vanishes, the degree of $F$ will be lower.
Let $\CC_{\hat F}$ and $\CC_{\hat G}$ be the zero sets of $\hat F$ and $\hat G$ in $\CP^2$.
Define the multiplicity of the solution to 
(\ref{eqbp1}) at the point $(\l,\mu)$ as the 
intersection multiplicity  of the curves $\CC_{\hat F}$ and $\CC_{\hat G}$
at the point  $(1:\l:\mu)$.

Equivalently
%, if (\ref{eqbp55}) is non-zero, 
one can define the multiplicity of the solution to 
(\ref{eqbp1}) at the point $(\l,\mu)$ as the 
intersection multiplicity  of the curves $\CC_F$ and $\CC_G$
at the point  $(\l,\mu)$.
%If (\ref{eqbp55}) is zero, replace $f$ and $g$ by 
%the Blaschke products $f \circ m_\gamma  $ and $g \circ m_\gamma $ instead,
%where $m_\gamma$ is the  M\"obius map
%\be
%\label{eqbp56}
%m_\gamma (z) \ := \ \frac{\gamma -z}{1 - \bar \gamma z}
%\ee
%and $\gamma$ is chosen as a non-critical point of $f$.

Define $F^\sim$ by
$$
F^\sim(z,w) \ := \ z^k w^\ell \overline{F ( 1/\bar z, 1/\bar w)},
$$
where $k$ is the highest power of $z$ in $F$ and $\ell$ is the highest power of $w$.
Then both $F$ and $G$ are symmetric, in the sense that 
$F^\sim$ equals $F$ and 
$G^\sim$ equals $G$. 
 
See \cite{ams06} for a discussion of polynomials that are symmetric with respect to the torus.
Moreover, since the functions $f$ and $g$ have modulus less than one in $\D$, equal to one on $\T := \partial \D$,
and greater than one on $\EE := \C \setminus \overline{\D}$,
both $\CC_F$ and $\CC_G$ are contained in $\D \times \D \, \cup \, \T \times \T \, \cup \, \EE \cup \EE$.

\bt\label{thmbp}
Let $f$ and $g$ be Blaschke products of degree $m$ and $n$ respectively, $m,n \geq 2$.
Using the notation (\ref{eqbp2} - \ref{eqbp5}), let
$N_{f,g}$ be
the number of points in $\CC_F \cap \CC_G$ in $\D \times \D$ plus half the number of points in 
$\T \times \T$.
Then $N_{f,g}$ 
is either infinite or $(m-1)(n-1)$.
Moreover, if  neither $\l$ nor $\mu$ is zero and  $(\l,\mu)$ is a solution to (\ref{eqbp1}), 
so is $(\ol,\om)$.
\et
\bp
Let us prove the last assertion first. It is an immediate consequence of the symmetry of $\CC_F$ and
$\CC_G$ with respect to $\T^2$. To see it directly, 
\beq
f(\ol) &\=&  \frac{p^\sim(\ol)}{p(\ol)} \\
&=&  \frac{\overline{p(\l)}}{\overline{p^\sim(\l)}} \\
&=& \frac{1}{\overline{f(\ol)}} .
\eeq
So if $f(\l) = f(\mu)$, then $f(\ol) = f(\om)$, and similarly for $g$.
\vs

We can assume without loss of generality that
(\ref{eqbp55}) is non-zero.
Indeed, otherwise we can choose some M\"obius map
\be
\label{eqbp56}
m_\gamma (z) \ := \ \frac{\gamma -z}{1 - \bar \gamma z}
\ee
and consider the Blaschke products $f \circ m_\gamma  $ and 
$g \circ m_\gamma $ instead, where $\gamma$ is chosen as a non-critical point of $f$.
%To see that some $\gamma$ can be found that renders (\ref{eqbp55}) non-zero, note that it
%can be rewritten as 
%\be
%\label{eqbp57}
%(-1)^{m-1} \left[ \sum \al_i - \sum \frac{1}{\bar \al_i} \right] \,
%\prod \bar \al_i .
%\ee
%If $\gamma$ is close to $1$, then all of the $m_\gamma (\al_i)$, which are the zeroes
%of $f \circ m_\gamma$, will lie in the orocycle
%$\D( 3/4, 1/4)$. Then 
%$$\Re \sum_{i=1}^m  m_\gamma (\al_i) \ < \ m $$
%and
%$$\Re \sum_{i=1}^m  \frac{1}{\overline{m_\gamma (\al_i)}} \ > \ m .$$
%Therefore (\ref{eqbp57}) cannot be $0$. 
%\vs
%
To count the intersections, first assume that $g(z) = z^n$, \ie that all the $\beta_i$ are $0$.
% and that no $\al_i = 0$.
For $\omega$ an $n^{\rm th}$ root of unity, $g(\l) = g(\omega \l)$.
For $\mu = \omega \l$ to also satisfy (\ref{eqbp1}), we need
\be
\label{eqbp6}
\frac{p(\l) p^\sim(\omega \l) - p(\omega \l) p^\sim(\l)}{\l - \omega \l} \= 0 .
\ee
The left-hand side of (\ref{eqbp6}) is a
polynomial of degree exactly $2m-2$. Indeed, the coefficient of $\l^{2m-2}$ is
$\omega^{m-1}$ times (\ref{eqbp55}), which we have assumed to be non-zero.
%\be
%\label{eqbp7}
%(-1)^{m-1} \omega^{m-1} \left[ \sum \al_i - \sum \frac{1}{\bar \al_i} \right] \,
%\prod \bar \al_i .
%\ee
%For generic $\alpha_i$, the coefficient (\ref{eqbp7}) is non-zero, and 
So for $n-1$ choices of
$\omega$, there are $2m-2$ different $\l$'s that solve (\ref{eqbp6}), giving
$2(m-1)(n-1)$ solutions to 
(\ref{eqbp1}) in $\C \times \C$. 

By the symmetry $F = F^\sim$, or by direct calculation, the constant term in 
$F(\l,\omega \l)$ is the complex conjugate of the highest order coefficient, and this
is non-zero by hypothesis. So none of the points $(\l,\omega\l)$ on $\CC_F$ have $\l = 0$.
Therefore each point with $\l \inn \D$ is matched by a point with $\l \inn \EE$, 
by the first part of the proof.
Therefore we have
$$
N_{f,z^n} \= \frac{1}{2} 2 (m-1) (n-1) \= (m-1)(n-1) .
$$

Now, move the zeroes of $g$ away from $0$ to the points $\beta_i$.
The function $G |_{\CC_F}$ 
is changing continuously, so the number of zeroes of $G$ in $\CC_F \cap \D^2$ can only
change either if a zero moves to $\CC_F \cap \T^2$, in which case its reflection also moves to 
$\CC_F \cap \T^2$ creating a double point there, or if a double point
 on $\CC_F \cap \T^2$ splits, and one moves in and the other moves
out.
In either case,  $N_{f,g}$ remains constant.

The only other possibility is that suddenly  $G |_{\CC_F}$ become identically zero
on some component.
This might happen at the end-point of the path, but can be avoided at all intermediate stages.
(Indeed, 
otherwise there would be a branch of an analytic function, $\omega(z) \not\equiv z$, such that
$f(z) \equiv f(\omega(z))$ and $g(z) \equiv g(\omega(z)$. Perturbing just one of $g$'s 
zeroes destroys this symmetry).

Therefore, at all points along the path, $N_{f,g} = (m-1)(n-1)$, and at the end-point
it either has this value or is infinite.
\ep

B\'ezout's theorem tells us that the number of intersection points in $\CP^2$
of $\hat F$ and $\hat G$ is precisely 
${\rm deg}(\hat F) \cdot {\rm deg}(\hat G)$. If one knew the multiplicities at infinity, one could
recover the number of intersections in $\C^2$. Unfortunately, this is hard to do; but one can use
Theorem~\ref{thmbp} to calculate the multiplicity at infinity.
\begin{cor}
\label{corbp1}
With notation as above, assume that $F$ and $G$ have degrees exactly $2m-2$ and $2n-2$, respectively,
and that $\CC_F \cap \CC_G$ is finite.
Let $r$ be the number of points $\l$ in $\D \setminus \{ 0 \}$ 
that satisfy $f(\l) = f(0)$ and $g(\l) = g (0)$. Let $s$ be the intersection multiplicity of 
$\CC_F$ and $\CC_G$ at $(0,0)$. Then the intersection multiplicity of the zero sets of $\hat F$ and
$\hat G$ at both $(0:1:0)$ and $(0:0:1)$ is 
exactly $$
(m-1)(n-1) + r + \frac{1}{2} s .
$$
\end{cor}
\bp
The number of solutions at points of the form $(1:z:w)$, by Theorem~\ref{thmbp}, is twice $N_{f,g}$ minus the
number of zeroes with either $z$ or $w$ equal to $0$. 
So the number at points of the form $(0:z:w)$ is, by  B\'ezout's theorem,
\be
\label{eqbp9}
4(m-1)(n-1) \ - \
2(m-1)(n-1) \ + \ 
2r \ + \ s .
\ee
At a point $(0:z:w)$, the homogeneous polynomial $\hat F$ is 
(\ref{eqbp55}) times $z^{m-1}w^{m-1}$, so the only roots at infinity occur when either
$z=0$ or $w=0$. By symmetry, the intersection multiplicities 
at $(0:1:0)$ and $(0:0:1)$ are the same, so each must be exactly half (\ref{eqbp9}).
\ep
\vs
We shall be  interested in proving that a pair of Blaschke products $(f,g)$
generates a specific
cofinite algebra (see Section~\ref{sece}). Recall that $\alg(f,g)$
denotes
the algebra generated by $f$ and $g$ (\ie the closure in $O(\D)$ of the set of polynomials in
$f$ and $g$).  
The codimension of   $\alg(f,g)$
is exactly half the number of points in 
$\CC_F \cap \CC_G \cap \D^2$ (because both $(\l,\mu)$ and $(\mu,\l)$ are counted).
The following theorem is the principal result in this section.
\bt\label{thmbp2}
Let $f$ and $g$ be Blaschke products of degree $m$ and $n$ respectively.
Suppose there are exactly $r < \infty$ unordered pairs of points on $\T$ that are not separated by $(f,g)$.
Then the codimension of $\alg(f,g)$ is 
$\frac{(m-1)(n-1) - r}{2} $.
\et
\bp
If the codimension were infinite, then there would be some distinguished variety on which both $F$ and $G$ vanished,
and hence there would be an infinite number of pairs of points in $\T \times \T$ that were identified by $(f,g)$.
When the codimension is finite, the result follows from Theorem~\ref{thmbp}.
\ep
\vs
Sometimes one can just prove that there are at least $r$ pairs of points on $\T$ that are identified
by $(f,g)$. To apply Theorem~\ref{thmbp2} to conclude that the codimension of $\alg(f,g)$ is at most
$\lfloor \frac{(m-1)(n-1) - r}{2} \rfloor$, one needs to know that 
$\alg(f,g)$ is cofinite. Here is a sufficient condition.
\begin{prop}\label{propbp1}
Let $f$ and $g$ be Blaschke products of degree $m$ and $n$ respectively, with $ 2 \leq m \leq n < \i$.
Using the notation (\ref{eqbp2} - \ref{eqbp5}), assume that $F$ is irreducible.
If $n$ is not a multiple of $m$, then $\alg(f,g)$ is cofinite. If $n$ is a multiple of $m$, then either
$\alg(f,g)$ is cofinite or $g$ is in $\alg(f)$.
\end{prop}
\bp
If $\CC_F  \cap \CC_G$ is infinite, then $G$ vanishes on all of $\CC_F$, and $F$ divides $G$ since $F$ is assumed
irreducible.
Therefore, whenever $f(\l) = f(\mu)$, we also have $g(\l)= g(\mu)$, so the function
$h := g \circ f^{-1}$ is well-defined and analytic, and $ g = h \circ f$ is in $\alg(f)$.
Since $f$ and $g$ are Blaschke products, so is $h$, and so the degree of $g$ is a multiple of the 
degree of $f$.
\ep
In practice, it may be hard to verify if $F$ is irreducible, so for later use we shall refine the proposition
in the case that $m$ is $2$ or $3$.
\begin{prop}\label{propbp2}
Let $f$ and $g$ be Blaschke products of degree $m$ and $n$ respectively, with $ 2 \leq m \leq n < \i$.
Assume that $m$ is either $2$ or $3$.
If $n$ is not a multiple of $m$, then $\alg(f,g)$ is cofinite. If $n$ is a multiple of $m$, then either
$\alg(f,g)$ is cofinite or $g$ is in $\alg(f)$.
\end{prop}
\bp
By composing $f$ and $g$ with some $m_\gamma$, as in the proof of Theorem~\ref{thmbp}, we can assume that
$F$ is of degree exactly $2(m-1)$ with the only term of that order being $z^{m-1}w^{m-1}$, and that
$G$ is of degree exactly $2(n-1)$ with the only term of that order being $z^{n-1}w^{n-1}$.

If $F$ is reducible, it cannot have a factor that is a polynomial in just $z$ (or just $w$), as this would force
$p$ to be constant. Moreover, it cannot have a linear factor, because if $a, b$ are non-zero and
$$
z^{m-1}w^{m-1} \ +\ {\rm lower\ order}\= ( az + bw + c)\, ({\rm stuff}) ,
$$
the second factor would have to have a term in either $z^{m-2} w^{m-1}$ or $z^{m-1} w^{m-2}$ to get a term
in $z^{m-1} w^{m-1}$. Either way, one would then also get a term in either $z^{m-2} w^m$ or $z^m w^{m-2}$, which
is not allowed.

So if $m=2$, we must have that $F$ is irreducible. If $m=3$, the only possible factoring is into a pair of quadratic factors.
If either factor has a term in say $z^2$, the other can only have powers of $w$, and this has been ruled out.
So the only possible factoring is of the form
$$
F(z,w) \= (a zw + b z + c w + d) \, ( a' zw + b' z + c'w + d') ,
$$
with $a a' \neq 0$.
This means $\CC_F$ decomposes into $\{(z,\phi_1(z))\} \cup \{ (z,\phi_2(z))\}$ for  linear fractional
transformations $\phi_1$ and $\phi_2$. As
$$
f(z) \= f(\phi_1 (z)) \= f(\phi_2 (z)) ,
$$
and these are the only points where this value is attained, 
%By symmetry about the diagonal, either $\phi_1 = \phi_1^{-1}$ and
%$\phi_2 = \phi_2^{-1}$ or 
we have $\phi_2 = \phi_1 \circ \phi_1$ and $\phi_1 \circ \phi_1 \circ \phi_1 = {\rm id}$.
Since $\CC_F$ is 
contained in $\D^2 \cup \T^2 \cup \EE^2$, it follows that $\phi_1$ and $\phi_2$ are disk automorphisms.

For definiteness, assume that $G$ vanishes on the sheet  $\{(z,\phi_1(z))\}$.
This means that
$$
g(z) \= g(\phi_1(z)) .
$$
%So $\phi$ permutes the zeroes of $g$. 
But then 
$$
g(z) \= g (\phi_1(z)) \= g(\phi_1\circ \phi_1 (z)),
$$
so $G$ actually vanishes on all of $\CC_F$.
Then we are back in the situation of the proof of Proposition~\ref{propbp1}, and we get
that $g$ is on $\alg (f)$.
\ep

\section{Petals}\label{secz}

A {\em petal} is a \gas that is holized by holomaps from $\D$, such as the set $V_1$
in Example~\ref{exa1} of Section~\ref{seca}. 
By Theorem~\ref{thmcc1}, two petals $V_1$ and $V_2$ are isomorphic if and only if there is a M\"obius
map $m: \D \to \D$ that pushes forward the first connection onto the second.

To determine if there is any holomorphic map from $V_1$ to $V_2$, again by Theorem~\ref{thmcc2} we can pass 
to their desingularizations, and ask if there is a holomorphic map $\psi$ from $ \D$ to $\D$
that pushes forward the connection $\Gamma_1$ corresponding to $V_1$ into the connection
$\Gamma_2$ corresponding to $V_2$.
The simplest case is when the only singularities are double points.
Then $\Gamma_1$ identifies the values at the points
$\alpha_{2i-1}$ and $\alpha_{2i}$, for $i = 1, \dots, n$, and 
$\Gamma_2$ identifies $\beta_{2j-1}$ and $\beta_{2j}$ for $ j = 1, \dots, m$, with $m \geq n$.
After reindexing the $\beta$'s, the problem becomes whether there is a map
$\psi : \D \to \D$ such that $\alpha_i$ goes to $\beta_i$ for every $ 1 \leq i \leq 2n$.
By Pick's theorem  \cite{pi16}, (see \eg
\cite{gar81}, \cite{foi-fra} or \cite{ampi} for an exposition), this can be done if and only if the matrix
$$
\left( \frac{1 - \overline{\beta_i} \beta_j}{1 - \overline{\alpha_i} \alpha_j}
\right)_{i,j=1}^{2n}
$$
is positive semi-definite. 

For more general connections, the question becomes whether there is a holomoprhic function 
$\psi : \D \to \D$ that satisfies a finite number of conditions on its values and its derivatives' values
on the support of $\Gamma_1$. Again, this is answerable using a modification of Pick's theorem
\cite{foi-fra} \cite{ampi}. 
In particular, if the condition is ever satisfied, it is always satisfied by a finite Blaschke product.
So we obtain:

\bprop
\label{propz1}
Let $V_1$ and $V_2$ be petals. If there exists any holomorphic map $\phi$ from
$V_1$ to $V_2$, then there exists a proper holomorphic map of constant finite valence.
\eprop

\section{Petals of codimension one}\label{secd}

In order to prove that certain holomaps actually realize 
particular algebras, the following lemma is useful (it can also be derived from Theorem~\ref{thm2}).

\bl
\label{lemcc2}
Suppose $\Gamma$ is an algebraic connection on $S$, and that $h$ is a holomap that holizes
the algebra $\Gamma^\perp$. Let $Dh$ denote the derivative of $h$.

\noindent
(i) If $\{ \alpha \} = \{ h^{-1}(h(\al)) \}$ and $Dh(\al) \neq 0$, then 
$\al \notin {\rm suppt}(\Gamma)$.

\noindent
(ii)
If $\{  h^{-1}(h(\al_1)) \} \= \{ \al_1, \dots, \al_m \}$
and there are functions $f_1, \dots, f_m$ in $A_h$ such that 
$Df_j(\al_i) = \delta_{ij}$, then the irreducible component of $\Gamma$
with support $\{ \al_1, \dots, \al_m \}$ is just
the set of equations
$$
f(\al_1) \= f(\al_2 ) \= \dots \= f(\al_m) .
$$
\el
\bp
(i) Suppose
\be
\label{eqcc2}
\Lambda(f) \= \sum_{i=1}^m \sum_{j=0}^{n_i} a_{ij}
f^{(j)}(\alpha_i) 
\ee
is a functional in $\Gamma$, with $\alpha = \al_1$.
Choose a polynomial $p$ so that
$p\circ h(\al) = 1,
\ D(p\circ h)(\al) = 1$ and $p\circ h (\al_i) = 0$ for all $i > 1$.
Then for every $k,\ell \geq 1$, the functions 
$[ (p\circ h )^k - 1 ]^\ell$ are in $A_h$, 
and the only way they can all be in the kernel of (\ref{eqcc2}) is
if all the $a_{ij}$'s are zero.

(ii) The functions $[ f_j^k - 1]^\ell$ give functions that vanish to high order at
each $\alpha_i$, for $i \neq j$, and the first  $\ell -1$ derivatives
vanish at $\al_j$. So if $\Lambda$ is as in (\ref{eqcc2}),
then every $a_{ij}$ must be zero if $j > 0$, and $\Lambda$ cannot actually depend on the values of the
derivatives.
\ep

The simplest petals that are not disks correspond to codimension one algebras.
There are two kinds. The first comes from a simple cusp. It can be described as the
algebra
\be\label{eqd1}
\{ f \inn O(\D) \ : \ f'(\alpha) = 0 \} ,
\ee
for some $\alpha \in \D$. As a \hacp, it can be 
holized by the map
\be\label{eqd2}
h(\zeta) \= ( m_\alpha (\zeta)^2, m_\alpha(\zeta)^3 ),
\ee
where $m_\alpha$ is the M\"obius map
$$
m_\alpha (\zeta) \ := \ \frac{\alpha - \zeta}{1 - \overline{\alpha} \zeta} .
$$
All of these algebras are isomorphic, and in fact each of the holizations 
in (\ref{eqd2}) actually holize the algebra into exactly the same set,
$$
\{ (z,w) \inn \C^2 \ : \ z^3 = w^2,\quad |z| < 1 \} .
$$
It follows therefore that if one takes the Neil parabola, that is the
curve $\{ z^3 = w^2 \}$,
and cuts it with any simple closed curve that contains $(0,0)$ in its interior,
the corresponding \hacs are isomorphic.
\vs
The second type of codimension one petal
is the single-crossing algebra defined, for each pair $(\alpha_1,\alpha_2)$ of distinct points in 
$\D$, as 
\be\label{eqd3}
A^1_{\alpha_1,\alpha_2} \ :=\ \{ f \inn O(\D) \ : \ f(\alpha_1) = f(\alpha_2) \} .
\ee
Two such algebras $A^1_{\alpha_1,\alpha_2}$ and $A^1_{\beta_1,\beta_2}$
are isomorphic if and only if there is an automorphism of the disk that takes the
pair $(\alpha_1,\alpha_2)$ to the pair $(\beta_1,\beta_2)$, and this in turn 
happens if and only if the hyperbolic
distance from $\alpha_1$ to $\alpha_2$ equals the 
 hyperbolic
distance from $\beta_1$ to $\beta_2$.

If one holizes the algebra $A^1_{\alpha_1,\alpha_2}$, for example by the map
$$
h(\zeta) \= ( m_{\alpha_1}(\z) m_{\alpha_2}(\z), \ \z m_{\alpha_1}(\z)
m_{\alpha_2}(\z)\ ),
$$
the hyperbolic geodesic connecting $\alpha_1$ to $\alpha_2$ gets mapped to the Kobayashi
geodesic that connects the multiple point to itself by going around the curve, and these
two geodesics have the same length. This length therefore forms a complete isomorphism
invariant for \gass with just a single crossing and no cusps, as in Example~\ref{exa1} of
Section~\ref{seca}.

\section{Petals of codimension two}\label{sece}

There are six different types of petals of codimension two. 

\subsection{Triple Points}

The first consists of those
functions that identify a triple $(\a_1,\a_2,\a_3)$ of distinct points:
\be\label{eqd4}
A^2_{\alpha_1,\alpha_2,\a_3} \ :=\ \{ f \inn O(\D) \ : \ f(\alpha_1) = f(\alpha_2) =
f(\a_3)\} .
\ee
The isomorphism problem for two such algebras $A^2_{\alpha_1,\alpha_2,\a_3} $ and
$A^2_{\beta_1,\beta_2,\beta_3}$ is easy: 
as explained in Section~\ref{secz},
%by Pick's theorem \cite{pi16}, (see \eg
%\cite{gar81} or \cite{ampi} for an exposition) 
there is a disk
automorphism that takes the first triple to the second if and only if one of
the $3$-by-$3$
matrices
$$
\left( \frac{1 - \overline{\sigma(\beta_i)} \sigma(\beta_j)}{1 - \overline{\alpha_i} \alpha_j}
\right)_{i,j=1}^3
$$
is rank one, where $\sigma$ ranges over the permutation group $S_3$. 

For the embedding dimension, note first that the local embedding dimension must exceed
two, because of the triple point. Indeed,
let $h$ be any holomap into $\C^2$ that sends $\alpha_1$, $\alpha_2$ and $\alpha_3$
to the same point, $\gamma$ say.
Then $\{ Dh (\alpha_i) \ :\  1\leq i \leq 3 \}$ is a set of three vectors in $\C^2$, so
must be linearly dependent. If $F$ is holomorphic in a neighborhood of $\gamma$,
then $\{ D(F \circ h) (\alpha_i) \}$ can be at most two dimensional, and therefore any
algebra holized by $h$ must have a linear relation on the derivatives at $\a_1, \a_2$
and $\a_3$. 

\bfig 
\resizebox{!}{2in}{\includegraphics {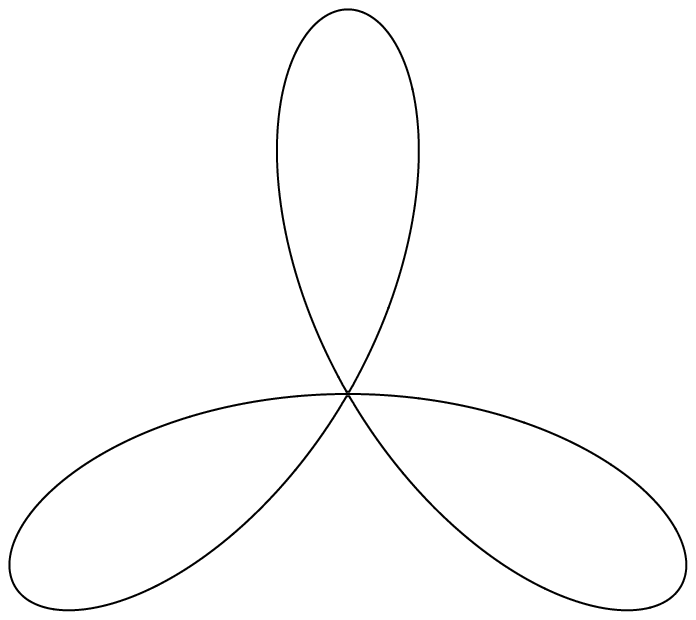}}
\caption{Petal with a triple point}
\label{pice1} \efig

A global embedding dimension of three can be attained. To see this, assume without loss
of generality that $\alpha_1 = 0$, and define $f(\z) = \z m_{\a_2}(\z) m_{\a_3}(\z)$.
Let $h$ be given by
$$
h(\z) \= (f(\z),\ \z f(\z),\, \z^2 f(\z) ) .
$$
Then $h$ is one-to-one except on $\{\a_1,\a_2,\a_3\}$. Its derivative never vanishes.
We just have to rule out the possibility of linear dependence between
the vectors $Dh(\alpha_i)$, and then we will be done by Lemma~\ref{lemcc2}.
 A calculation yields that the matrix whose rows
are $Dh(\alpha_i)$ is
$$
\left(
\begin{array}{ccc}
\a_2 \a_3 & 0 & 0 \\
\frac{- \a_2 m_{\a_3}(\a_2)}{1 - |\a_2|^2} &
\frac{- \a_2^2 m_{\a_3}(\a_2)}{1 - |\a_2|^2} &
\frac{- \a_2^3 m_{\a_3}(\a_2)}{1 - |\a_2|^2} \\
\frac{- \a_3 m_{\a_2}(\a_3)}{1 - |\a_3|^2} &
\frac{- \a_3^2 m_{\a_2}(\a_3)}{1 - |\a_3|^2} &
\frac{- \a_3^3 m_{\a_2}(\a_3)}{1 - |\a_3|^2} 
\end{array} \right).
$$
This has non-zero determinant, so has full rank.

\subsection{Cusps of Order 2}

Consider now a connection supported on one point, which we may as well take
to be the origin. A cofinite algebra of codimension 2 coming from such a connection
must be contained in an algebra of codimension 1 by Theorem~\ref{thm2}, 
so one of the linear functionals must be
$$
f \ \mapsto \ f'(0) .
$$
The other must be of the form
\be
\label{eqe2}
\Lambda \, : \, f \mapsto \ \sum_{j=2}^n a_j f^{(j)} (0) .
\ee
If $a_2 = 0$, then $z^2$ is in the algebra, and all the non-zero $a_j$ must have
$j$ odd. 
If only one $a_j \neq 0$, it must be $a_3$, and we get the algebra
$$
\{ f \inn O(D) \ : \ f'(0) = f'''(0) = 0 \}.
$$
This algebra is holized by the functions $(z^2,z^5)$.

Otherwise, let $j_0$ be the index of the first non-zero coefficient, and $n \geq j_0 +2$
be the index of the last.
The algebra contains some polynomial 
of the form
$$
p(z) \= z^{j_0} + c_{j_0 +2} z^{j_0 +2} + \dots + c_n z^n
$$
Subtracting appropriate linear multiples of $z^2 p(z), z^4 p(z), ...$
from $p$, one gets that $z^{j_0}$ is in the algebra, and so
$\Lambda$ has to be the zero functional, a contradiction.

So assume $a_2 \neq 0$. Suppose $k$ is the smallest integer larger than 1 such that $z^k$ is not in the algebra.
If $k > 2$, then
$$
p(z) \ := \ z^2 - \frac{2\, a_2}{k!\, a_k} z^k
$$
is in the algebra, 
and so is $p^2$, and hence so is 
$z^4$. Therefore $k$ is either two or three.

Putting all this together, we conclude that the algebra must be of the form
\be
\label{eqe3}
B^2_{0,c} \ := \ \{ f \in O(\D) \ : \ f'(0) =0,\ f''(0) + c f'''(0) = 0 \}
\ee
for some $c$ in $\C$. Two such algebras are isomorphic iff $|c|$ is the same.

If $c = 0$, then the algebra is not locally embeddable in $\C^2$. Indeed, however the generators are chosen, 
there must be a linear relationship between $f'''(0), f^{(iv)}(0)$ and $f^{(v)}(0)$.

If $c \neq 0$, the algebra can be holized in $\C^2$, 
%- see Section~\ref{sech}. 
as we shall show in a subsequent paper \cite{amcusps}.
For $|c|$ large enough,
the holization has a particularly simple form.
Let $f(z) = z^4$ and $g(z) = z^2 m_\alpha(z)$, where $\alpha$ is chosen so that
\be
\label{eqe5}
g(1) = g(i). 
\ee
A calculation shows that
(\ref{eqe5}) holds iff $\alpha$ lies on the circle centered at $1+i$ and radius $1$.
Another calculation shows that $g$ lies in $B^2_{0,c}$ for 
\be
\label{eqe6}
c \= \frac{\alpha}{3(1-|\alpha|^2)} .
\ee
So for this value of $c$, both $f$ and $g$ lie in  $B^2_{0,c}$, and they identify the points $1$ and $i$.
Therefore, by Theorem~\ref{thmbp2} and Proposition~\ref{propbp2}, they generate the algebra  $B^2_{0,c}$.
As $\alpha$ moves along the admissible circular arc, we get that every algebra 
$B^2_{0,c}$ is holized by $z^4$ and a suitable Blascke product of degree $3$
provided that
$$
|c| \ \geq \ \frac{2-\sqrt{2}}{3(4\sqrt{2}-5)} \= 0.297\cdots .
$$

\subsection{Other codimension 2 petals}

There are four other types of codimension 2 petals: two crossings, two cusps, and a crossing and a cusp
which may or may not be at one of the points on the crossing.
The generic such petal is two crossings, and for four distinct points $\al_1,\al_2,\beta_1,\beta_2$
we shall let
$$
C^2_{\al_1,\al_2;\beta_1,\beta_2} \ := \ \{ f \inn O(\D) \ : \ f(\al_1) = f(\al_2),\ f(\beta_1) = f(\beta_2) \}.
$$

If one wishes to map the closed disk $\overline{\D}$ into $\C^2$ with a pair of Blaschke products so that
there are precisely two crossings, then the points $\al_1,\al_2,\beta_1,\beta_2$ must satisfy a polynomial 
equation.
\bt
\label{thme5}
The algebra $C^2_{\al_1,\al_2;\beta_1,\beta_2} $
can be holized by a pair of finite Blaschke products $(f,g)$
that are also  one-to-one on $\T$ if and only if
\be
\label{eqe8}
\frac{(\al_1 - \beta_1)(\al_2 - \beta_1)} {( 1 - \overline{\al_1}\beta_1 ) ( 1- \overline{\al_2} \beta_1)}
\=
\frac{(\al_1 - \beta_2)(\al_2 - \beta_2)}{ ( 1 - \overline{\al_1}\beta_2 ) ( 1- \overline{\al_2} \beta_2)}
.\ee
\et
\bp
Let $f = m_{\al_1} m_{\al_2}$. Condition (\ref{eqe8}) is the 
requirement that $f(\beta_1) = f(\beta_2)$.
If this holds, then $f$ is a Blaschke product of degree two that lies in the algebra.
Let $g$ be any Blaschke product of degree 5 in the algebra, for example
$$
g(z) = z m_{\al_1}(z) m_{\al_2}(z) m_{\beta_1}(z)m_{\beta_2}(z).$$
Then by Theorem~\ref{thmbp2} and Proposition~\ref{propbp2}, the pair $f$ and $g$ holize the algebra
and do not identify any points on $\T$.
\vs
Conversely, suppose $f$ and $g$ are Blaschke products that holize the algebra and are one-to-one on 
$\T$. By Theorem~\ref{thmbp2}, they either have degrees $2$ and $5$ or degrees $3$ and $3$.
In the first case, if $f$ is in the algebra and has degree $2$, let $f(\alpha_1) = \gamma$.
Then $m_\gamma \circ f$ is a constant times $ m_{\al_1} m_{\al_2}$, and has the same values
at $\beta_1$ and $\beta_2$, so (\ref{eqe8}) holds.

The second case cannot hold. Indeed, suppose $f$ and $g$ are both Blaschke products of degree 3 in the algebra.
After composing with the M\"obius transformations that send their values at $\alpha_1$ to $0$, we 
can take
\beq
f(z) &\=& m_{\al_1}(z) m_{\al_2}(z) m_\gamma(z) \\
g(z) &\=& m_{\al_1}(z) m_{\al_2}(z) m_\delta(z) .
\eeq
Let $ m_{\al_1}(\beta_1) m_{\al_2}(\beta_1) =  C_1$ and 
$ m_{\al_1}(\beta_2) m_{\al_2}(\beta_2) =  C_2$.
Then
\se\att
\begin{eqnarray}
\label{eqe9}
C_1 m_\gamma(\beta_1) &\=& C_2 m_\gamma(\beta_2) \\
C_1 m_\delta(\beta_1) &\=& C_2 m_\delta(\beta_2) .
\nonumber
\end{eqnarray}
But there is a unique $\gamma$ satisfying (\ref{eqe9}), as can be seen either by direct calculation,
or observing that $m_\gamma$ is the unique solution to the extremal problem of finding a 
holomorphic map from $\D$ to $\D$
that sends $\beta_1$ to $r$ and $\beta_2$ to $\frac{C_2}{C_1} r$, with $r$ as large as possible.
So $g$ would be a function of $f$, and the pair would not holize the algebra.
\ep
\vs
The algebra that has a single crossing and a cusp at the crossing, such as
$$
\{ f \inn O(\D) \ : \ f(0) = f(\al),\ f'(0) = 0 \} 
$$
can not be locally holized in two dimensions. The other three can be.
\vs
We do not know the answer to either of the following questions, even in the case
of codimension 2. 
\begin{question}
If a petal can be locally holized in dimension $2$, can it be globally
holized by a pair of functions?
\end{question}
\vs
\begin{question}
If a petal can be holized by a pair of functions, can it be holized by a pair
of Blaschke products?
\end{question}

%The simplest petals are those that
%come from a connection of order $0$ (\ie the linear functionals depend
%only on the function values, not their derivatives). One then gets an
%algebra of the form
%\be
%\label{eqd1}
%\{ f \inn O(\D) \ : \ f(\alpha_{ij})\ =\ f(\alpha_{ik}), \quad 1 \leq i
%\leq N,\ 1 \leq j,k \leq N_i \} ,
%\ee
%for $N$ sets $\{\alpha_{i1},\dots,\alpha_{i N_i} \}$ of points in $\D$,
%with all the points distinct.

%\section{Cusps}\label{sech}

\section{ Model theory}\label{seci}

In this section, we propose to study what happens to the Hardy space theory
if the disk is replaced by a nice \hacp.
When there are no singularities, the Hardy space theory for bordered Riemann surfaces has been 
developed by J.~Ball and V.~Vinnikov \cite{bv00}.
Throughout this section, we shall take $\CC$ to be the algebraic set
$$
\CC \= \{ (z,w) \inn \C^2 \ : \ z^2 = w^2 \} ,
$$
and $V$ to be $\CC \cap \D^2$.
The set $V$ consists of two disks:
\beq
\D^+ &\ := \ & \{ (\z,\z) \ : \ \z \inn \D \} \\
\D^-  &\ := \ & \{ (\z,-\z) \ : \ \z \inn \D \} .
\eeq
If we let $0^+$ and $0^-$ denote the centers of $\Dp$ and $\Dm$, respectively, we can think
of $O(V)$ as the set of holomorphic functions $f$
on two disjoint disks with the connection $f(\zp) = f(\zm)$.

Let $\sp$ and $\sm$ denote Lebesgue measure on $\partial \Dp$ and $\partial \Dm$ respectively.
If one replaces the defining equation for $\CC$ by $z^2 = m_{\al_1} (w) m_{\al_2} (w) $
for distinct points $\alpha_1$ and $\alpha_2$, then $\CC \cap \D^2$ is an annulus.
In the limiting case of $V$, it  behaves like an annulus function theoretically,
in the sense that the obstruction to solving certain problems is give by one real parameter.

\bprop
\label{propia}
(i) Let $u$ be a real-valued function 
on $\partial V$. Then $u$ is the real part of a holomorphic function on $V$ iff
$$
\int_{\partial \Dp} u d \sp \= 
\int_{\partial \Dm} u d \sm .
$$

(ii) The points $\{ (\al_1,\al_1) ,\dots, (\al_m,\al_m), 
(\beta_1,-\beta_1), \dots, (\beta_n,-\beta_n) \}$ 
in $V \setminus \{(0,0)\}$ form the zero set of a rational inner function
iff
$$
\prod_{i=1}^m |\al_i| \= 
\prod_{j=1}^n |\beta_j| .
$$
\eprop
The proof is straightforward, and we omit it.

\vs
Let $T = (T_1,T_2)$ be a pair of commuting operators on a Hilbert space.
We say that $T$ has $V$ as a {\em spectral set}
if, for every polynomial $p$ in two variables,
\be
\label{eqi1}
\| p(T_1, T_2) \| \ \leq \ \| p \|_V .
\ee
Von Neumann's inequality \cite{vonN51} says that a single operator has the disk as a spectral set
iff it is a contraction. The analogous theorem for $V$ has a parameter in it.
Note first that if $T$ has $V$ as a spectral set, then in particular any polynomial that vanishes on
$V$ must vanish on $T$, so 
\be
\label{eqi2}
T_1^2 \= T_2^2. 
\ee
If (\ref{eqi2}) holds, we shall write $T \precc V$.
\bt
\label{thmia}
Suppose  $T \precc V$. Then $T$ has $V$ as a spectral set iff
\be
\label{eqi3}
\| T_1 + T_2 + e^{i\theta} (T_1 - T_2) \| \ \leq \ 2 \qquad \forall\ \theta \inn \R .
\ee
\et
\bp
By taking the Cayley transforms 
of the functions in (\ref{eqi1}), proving that $V$ is a spectral set is the same as
proving that whenever a holomorphic function has positive real part on $V$, then it has positive real part on
$T$. To analyze this condition, it is sufficient to look at the extreme points of the maps
$\phi$ from $V$ to the right-half plane, normalized to take the origin to $1$.
These extreme points are the functions
\beq
\phi |_{\Dp} \, : \, (\z,\z) &\ \mapsto \ & \frac{\tau^+ + \z}{\tau^+ - \z} \\
\phi |_{\Dm} \, : \, (\z,-\z) &\ \mapsto \ & \frac{\tau^- + \z}{\tau^- - \z} ,
\eeq
where $\tau^+$ and $\tau^-$ are unimodular.

Now, 
\be
\label{eqi6}
\phi(T) \= \phi |_{\Dp} (\frac{1}{2}(T_1 + T_2)) \ + \
\phi |_{\Dm} (\frac{1}{2}(T_1 - T_2)) - \phi(0) .
\ee
Then (\ref{eqi6}) has positive real part iff
(\ref{eqi3}) holds with $e^{i\theta} = \tau^+/\tau^-$.
\ep

\rem
Condition (\ref{eqi3}) can be rephrased in terms of the numerical radius of a single
operator, because of the following lemma whose proof we omit (if the first operator
in (\ref{eqi5}) is not invertible, take an appropriate limit).
\bl
\label{lemia}
Operators $A$ and $B$ satisfy
\be
\label{eqi4}
\| A  + e^{i\theta} B \| \ \leq \ 1 \qquad \forall\ \theta \inn \R 
\ee
iff the numerical radius of
\be
\label{eqi5}
( 1 - A^\ast A - B^\ast B)^{-1/2}  \, A^\ast B \,
( 1 - A^\ast A - B^\ast B)^{-1/2}
\ee
is less than or equal to $\frac{1}{2}$.
\el

\vs
The von Neumann-Wold theorem states that all isometries are direct sums of unitaries and copies
of the forward shift. What do pairs of isometries that live on $V$ look like?

\bt
\label{thmi2}
Let $T = (T_1,T_2)$ be a pair of commuting isometries on a Hilbert space $\h$ satisfying $T_1^2 = T_2^2$.
Then there is a decomposition 
$\h \= \M_{+} \oplus \M_{-} \oplus \K$, isometries $W_{\pm}$ on 
$\M_{\pm}$, and a pair of operators $E_{\pm} : \K \to \M_{\pm}$ satisfying
$\Ep^\ast \Ep + \Em^\ast \Em = I$, $\Wp^\ast \Ep = 0$, and $\Wm^\ast \Em = 0$,
so that
\se\att
\begin{eqnarray}
\nonumber
T_1 &\ \cong \ & 
\bordermatrix{&\Mp  &\Mn& \K\cr
\Mp &\Wp & 0&\Ep\cr
\Mn&0&\Wm &\Em\cr
\K&0&0&0} \\
\nonumber&&\\
\nonumber&&\\
T_2 &\ \cong \ & 
\bordermatrix{&\Mp  &\Mn& \K\cr
\Mp &\Wp & 0&\Ep\cr
\Mn&0&-\Wm &-\Em\cr
\K&0&0&0} .
\label{eqi7}
\end{eqnarray}
\et
\bp
Let $\Mp$ be the range of $T_1 + T_2$, and $\Mn$ be the range of $T_1 -T_2$.

Claim: $\Mp$ and $\Mn$ are orthogonal.

Proof of Claim:
Suppose first that $T_1$ and $T_2$ are unitary. 
\beq
(T_1 - T_2)^\ast (T_1 + T_2) T_1 &\=&
(T_1^\ast T_2 - T_2^\ast T_1) T_1 \\
&=& T_1^\ast T_1 T_2 - T_2^\ast T_2 T_2 \\
&=& 0 .
\eeq
As $T_1$ is unitary, it follows that 
$$
(T_1 - T_2)^\ast (T_1 + T_2)  \= 0 ,
$$
so the range of $T_1 - T_2$ is orthogonal to the range of $T_1 + T_2$.

For general isometries $T_1$ and $T_2$ as in the statement of the theorem, there
are commuting  unitaries $U_1$ and $U_2$ that extend them and satisfy $U_1^2 = U_2^2$
(as can be seen, for example, from the Sz.-Nagy-Foia\c{s} model for commuting isometries \cite{szn-foi}).
From the previous argument, the range of $U_1 - U_2$ is orthogonal to the range of $U_1 + U_2$, so
the range of $T_1 - T_2$ is orthogonal to the range of $T_1 + T_2$.
\vs
Let $\K := \h \ominus (\Mp \oplus \Mn)$.
Then $\K$ is orthogonal to the range of $T_1$ and $T_2$, so the third rows of
(\ref{eqi7}) will be zero. Using the fact that $ T \precc V$  and that 
$T_1$ and $T_2$ are isometries, one gets the rest of (\ref{eqi7}).
\ep

\section{Neil parabola}\label{secnp}

We shall give another example of how the function theory interacts with the geometry of
a \hacp.
Let
$$
V_2 \= \{ (z,w) \inn \D^2 \ : \ z^3 = w^2 \} .
$$
In a recent paper \cite{kn07a}, G.~Knese succeeded in finding an exact formula
for the Carath\'eodry metric on $V_2$; this was the first example of an explicit calculation
of the Carath\'eodory metric on a singular \hacp.

What are the extreme points of the set of holomorphic functions on $V_2$
with positive real part, normalized to take $(0,0)$ to $1$?
As $V_2$ is a petal, holized by the map
$$
h_2 \ : \ z \ \mapsto \ (z^2, z^3 ) ,
$$
the question is equivalent to asking for the extreme points of the 
set $\P$ of functions on the disk that have positive real part, send 0 to 1, and whose derivative
vanishes at the origin.
Every function $\phi$ in $\P$ has a Herglotz representation as
$$
\phi(z) \= \int \frac{e^{i\theta} + z}{e^{i\theta} - z} d \mu(\theta)
$$
for some probability measure $\mu$ on the circle. Being in $\P$
means in addition that 
\be
\label{eqnp1}
\int e^{-i\theta} d \mu(\theta) \= 0 .
\ee
Absent (\ref{eqnp1}), the extreme points just correspond to measures that have a single atom.
With the condition, there are two sorts of extreme point:
the first is two atoms of mass $1/2$ that are diametrically opposite.
The second is three atoms that form the vertices of a triangle with $0$ in the interior, and
such that their masses sum to $1$ and (\ref{eqnp1}) is satisfied.
One never needs more than three atoms, since by a theorem of Carath\'eodory if a point in $\R^d$
is in a convex set, it is in the convex hull of $d+1$ extreme points
(see \eg \cite{ga00} for a proof).

After taking the Cayley transform, we see that the minimal inner functions on $V_2$ either have two 
zeroes
(which must be $(\al^2,\al^3)$ and $(\al^2,-\al^3)$),
or three.

\section{Conclusion}\label{secg}

The purpose of this paper is to show how studying \hacs gives rise to a host of problems
in classical function theory, in Riemann surfaces, in function algebras and in operator theory.

A cofinite algebra has local restrictions on the minimal dimension of the space
in which it can be holized, as in Example~\ref{exa2}.
There are also global restrictions -- the desingularizing Riemann surface $S$ will not in
general be embeddable in $\C^1$. It is still an open question whether
all finite Riemann surfaces are embeddable in $\C^2$ \cite[p. 347]{bir66}. 
Let the global embedding dimension of a cofinite algebra $A$ 
be the smallest $n$ such that $A$ can be holized in $\cn$.
Let the local embedding dimension at any point $\lambda$ in the maximal ideal
space be the smallest $n$ such that some neighborhood of $\lambda$
can be embedded in $\cn$.
Perhaps the most significant open problem arising from our approach is the following. 
\begin{question}
When is the 
global embedding dimension of a cofinite algebra on $S$
equal to the maximum of the local embedding
dimensions and the global embedding dimension of $S$?
\end{question}

\bibliography{references}

\end{document}